\documentclass[10pt]{article}
\usepackage{amsmath}
\usepackage{latexsym}
\usepackage{amssymb}
\usepackage{amsfonts}
\usepackage{amsthm}
\title{DIFFERENTIAL HOMOLOGICAL ALGEBRA \\ AND GENERAL RELATIVITY}
\author{J.-F. Pommaret \\ CERMICS, Ecole des Ponts ParisTech  \\  jean-francois.pommaret@wanadoo.fr  \\
 http://cermics.enpc.fr/$\sim$pommaret/home.html }
\date{  }
\begin{document}
\maketitle

\noindent
{\bf ABSTRACT}:  \\
In 1916, F.S. Macaulay developed specific localization techniques for dealing with "unmixed polynomial ideals" in commutative algebra, transforming them into what he called "inverse systems" of partial differential equations. In 1970, D.C. Spencer and coworkers studied the formal theory of such systems, using methods of homological algebra that were giving rise to "differential homological algebra", replacing unmixed polynomial ideals by "pure differential modules". The use of "extension modules" and "differential double duality" is essential for such a purpose. In particular, 0-pure differential modules are torsion-free and admit an "absolute parametrization" by means of arbitrary potential like functions. In 2012, we have been able to extend this result to arbitrary pure modules, introducing a "relative parametrization" where the potentials should satisfy compatible "differential constraints". We recently discovered that General Relativity is just a way to parametrize the Cauchy stress equations by means of the formal adjoint of the Ricci operator in order to obtain a "minimum parametrization" by adding sufficiently many compatible differential constraints, exactly like the Lorenz condition in electromagnetism. These unusual purely mathematical results are illustrated by many explicit examples and even strengthen the comments we recently provided on the mathematical foundations of General Relativity and Gauge Theory.  \\

\vspace{2cm}
purity filtration; inverse system; involution; control theory; electromagnetism; General relativity.  \\

\newpage

\noindent
{\bf 1) INTRODUCTION}:  \\
Te main purpose of this paper is to prove how apparently totally abstract mathematical tools, ranging among the most difficult ones of differential geometry and homological algebra, can also become useful and enlight many engineering or physical concepts (See the review Zbl 1079.93001 for the only application to control theory). \\
In the second section, we first sketch and then recall the main (difficult) mathematical results on {\it differential extension modules} and {\it differential double duality} that are absolutely needed in order to understand the {\it purity} concept and, in particular, the so-called {\it purity filtration} of a differential module ([1],[2],[24],[29]). We also explain the unexpected link existing between {\it involutivity} and {\it purity} allowing to exhibit a {\it relative parametrization} of a pure differential module, even defined by a system of linear PD equations with coefficients in a non-constant differential field $K$. It is important to notice that the {\it reduced Spencer form} which is used for such a purpose generalizes the {\it Kalman form} existing for an OD classical control system and we shall illustrate this fact. \\
The third section will present for the first time a few explicit motivating academic examples in order to illustrate the above mathematical results, in particular the unexpected striking situations met in the study of contact and unimodular contact structures.  \\
In the fourth section, we finally provide examples of applications, studying the mathematical foundations of OD/PD control theory ([24],[25]), electromagnetism (EM) ([31],[37]) and general relativity (GR) ([30],[33],[35]). 
Most of these examples can be now used as test examples for certain computer algebra packages recently developped for such a purpose ([43-44]).  \\

\noindent
{\bf 2) MATEMATICAL TOOLS}:  \\

Let $D=K[d_1,...,d_n]=K[d]$ be the ring of differential operators with coefficients in a differential field $K$ of characteristic zero, that is such that $\mathbb{Q}\subset K$, with $n$ commuting derivations ${\partial}_1,...,{\partial}_n$ and commutation relations $d_ia=ad_i+{\partial}_ia,\forall a\in K$. If $y^1,...,y^m$ are $m$ differential indeterminates, we may identify $Dy^1+...+Dy^m=Dy$ with $D^m$ and consider the finitely presented left differential module $M={}_DM$ with presentation $D^p\rightarrow D^m\rightarrow M \rightarrow 0$ determined by a given linear multidimensional system with $n$ independent variables, $m$ unknowns and $p$ equations. Applying the functor $hom_D(\bullet,D)$, we get the exact sequence $0\rightarrow hom_D(M,D)\rightarrow D^m\rightarrow D^p \longrightarrow N_D \longrightarrow 0$ of {\it right differential modules} that can be transformed by a side-changing functor to an exact sequence of finitely generated {\it left differential modules}. This new presentation corresponds to the {\it formal adjoint} $ad({\cal{D}})$ of the linear differential operator $\cal{D}$ determined by the initial presentation but now with $p$ unknowns and $m$ equations, obtaining therefore a new finitely generated {\it left differential module} $N={}_DN$ and we may consider $hom_D(M,D)$ as the {\it module of equations} of the {\it compatibility conditions} (CC) of $ad({\cal{D}})$, a result not evident at first sight (See [24],[38]). Using now a maximum free submodule $0 \longrightarrow D^l \longrightarrow hom_D(M,D)$ and repeating this standard procedure while using the well known fact that $ad(ad({\cal{D}}))={\cal{D}}$, we obtain therefore an embedding $0\rightarrow hom_D(hom_D(M,D),D)\rightarrow D^l$ of left differential modules for a certain integer $1\leq l<m$ because $K$ is a field and thus $D$ is a noetherian bimodule over itself, a result leading to $l=rk_D(hom_D(M,D))=rk_D(M)< m$ as in ([22], p 341,[23], [25] p 179)(See section 3 for the definition of the {\it differential rank} $rk_D$). Now, the kernel of the map $\epsilon:M\rightarrow hom_D(hom_D(M,D),D):m\rightarrow \epsilon(m)(f)= f(m),\forall f\in hom_D(M,D)$ is the torsion submodule $t(M)\subseteq M$ and $\epsilon$ is injective if and only if $M$ is torsion-free, that is $t(M)=0$. In that case, we obtain by composition an embedding $0\rightarrow M \rightarrow D^l$ of $M$ into a free module (that can also be obtained by localization if we introduce the ring of fractions $S^{-1}D=DS^{-1}$ when $S=D-\{0\}$). This result is quite important for applications as it provides a (minimal) parametrization of the linear differential operator $\cal{D}$ and amounts to the controllability of a classical control system when $n=1$ ([12],[24]). This parametrization will be called an "{\it absolute parametrization} " as it only involves arbitrary "{\it potential-like} " functions (See [23],[26],[29],[30],[33],[36],[39],[41] for more details and examples, in particular that of Einstein equations).   \\

The purpose of this paper is to extend such a result to a much more general situation, that is {\it when} $M$ {\it is not torsion-free}, by using unexpected results first found by F.S. Macaulay in 1916 ([15]) through his study of "{\it inverse systems} " for "{\it unmixed polynomial ideals} ".  \\

Introducing $t_r(M)=\{m\in M\mid cd(Dm)>r\}$ where the codimension of $Dm$ is $n$ minus the dimension of the characteristic variety determined by $m$ in the corresponding system for one unknown, we may define the {\it purity filtration} as in ([1],[24],[29]):\\
    \[     0=t_n(M) \subseteq t_{n-1}(M) \subseteq  ... \subseteq t_1(M) \subseteq t_0(M)=t(M)\subseteq M  \]
The module $M$ is said to be $r$-{\it pure} if $t_r(M)=0, t_{r-1}(M)=M$ or, equivalently, if $cd(M)=cd(N)=r, \forall N\subset M$ and a torsion-free module is a 0-pure module. Moreover, when $K=k=cst(K)$ is a field of constants and $m=1$, a pure module is {\it unmixed} in the sense of Macaulay, that is defined by an ideal having an equidimensional primary decomposition. \\

\noindent
{\bf Example 2.1} :  As an elementary example with $K=k=\mathbb{Q}, m=1,n=2, p=2$, the differential module defined by $d_{22}y=0,d_{12}y=0$ is not pure because $z'=d_2y$ satisfies $d_2z'=0,d_1z'=0$ while $z"=d_1y$ only satisfies $d_2z"=0$ and $(({\chi}_2)^2,{\chi}_1{\chi}_2)=({\chi}_2)\cap ({\chi}_1,{\chi}_2)^2$. We obtain therefore the purity filtration $ 0 = t_2(M) \subset  t_1(M) \subset t_0(M)=t(M)=M $ with strict inclusions as $ 0\neq z' \in t_1(M) $ while $z" \in t_0(M)$ but $z" \notin t_1(M) $.   \\

From the few (difficult) references ([[1],[2],[3],[10],[13],[14],[16],[17],[18],[24],[25],[29],[45],[46]) dealing with the extension modules $ext^r(M)=ext^r_D(M,D)$ and purity in the framework of algebraic analysis, it is known that $M$ is $r$-pure if and only if there is an embedding $0 \rightarrow M \rightarrow ext^r_D(ext^r_D(M,D),D)$. Indeed, the case $r=0$ 
is exactly the one already considered because $ext^0_D(M,D)=hom_D(M,D)$ and the ker/coker exact sequence ([25],[29]):\\
\[   0 \longrightarrow ext^1(N) \longrightarrow M \longrightarrow ext^0(ext^0(M)) \longrightarrow ext^2(N) \longrightarrow 0  \]
allows to test the torsion-free property of $M$ in actual practice by using the double-duality formula $t(M)=ext^1(N)$ as in ([24],[25]).

Independently of the previous results, the following procedure, {\it where one may have to change linearly the independent variables if necessary}, is the heart towards the next effective definition of involution. It is intrinsic even though it must be checked in a particular coordinate system called $\delta$-{\it regular} ([19],[20],[32]) and is quite simple for first order systems without zero order equations.    \\

\noindent
$\bullet$ {\it Equations of class} $n$: Solve the maximum number ${\beta}^n_q$ of equations with respect to the jets of order $q$ and class $n$. Then call $(x^1,...,x^n)$ {\it multiplicative variables}.\\
\noindent
$\bullet$ {\it Equations of class} $i\geq 1$: Solve the maximum number ${\beta}^i_q$ of {\it remaining} equations with respect to the jets of order $q$ and class $i$. Then call $(x^1,...,x^i)$ {\it multiplicative variables} and $(x^{i+1},...,x^n)$ {\it non-multiplicative variables}.\\
\noindent
$\bullet$ {\it Remaining equations equations of order} $\leq q-1$: Call $(x^1,...,x^n)$ {\it non-multiplicative variables}.\\

\noindent
In actual practice, we shall use a {\it Janet tabular} where the multiplicative "variables" are represented by their index in upper left position while the non-multiplicative variables are represented by dots in lower right position ([11],[19],[24]) (Compare to ([47]).  \\

\noindent
{\bf DEFINITION 2.2}: A system of PD equations is said to be {\it involutive} if its first prolongation can be achieved by prolonging its equations only with respect to the corresponding multiplicative variables. In that case, we may introduce the {\it Cartan characters} ${\alpha}^i_q=m\frac{(q+n-i-1)!}{(q-1)!((n-i)!}-{\beta}^i_q$ for $i=1, ..., n$ and we have $dim(g_q)=\sum {\alpha}_q={\alpha}^1_q+...+{\alpha}^n_q$ and $dim(g_{q+1})=\sum i{\alpha}^i_q=1{\alpha}^1_q + ... +n{\alpha}^n_q$. Moreover, one can exhibit the {\it Hilbert polynomial} $dim(R_{q+r})$ in $r$ with leading term $(\alpha/d!)r^d$ with $d\leq n$ when $\alpha$ is the smallest non-zero character in the case of an involutive symbol. Such a prolongation allows to compute {\it in a unique way} the principal ($pri$) jets from the parametric ($par$) other ones. This definition may also be applied to nonlinear systems as well.  \\

\noindent
{\bf REMARK 2.3}: For an involutive system with $\beta ={\beta}^n_q< m$, then $(y^{\beta +1},...,y^m)$ can be given arbitrarily and may constitute the {\it input} variables in control theory, though it is not necessary to make such a choice. {\it In this case}, the intrinsic number $\alpha={\alpha}^n_q=m-\beta> 0$ is called the $n$-{\it character} and is the system counterpart of the so-called "{\it differential transcendence degree}" in differential algebra and the "{\it rank}" in module theory. As we shall see in the next Section, {\it the  smallest non-zero character and the number of zero characters are intrinsic numbers that can most easily be known by bringing the system to involution} and we have 
${\alpha}^1_q\geq ... \geq {\alpha}^n_q\geq 0$. \\

In the situation of the last remark, the following procedure will generalize for PD control systems the well known first order Kalman form of OD control systems where the derivatives of the input do not appear ([24], VI, Remark 1.14, p 802). For this, {\it we just need to modify the Spencer form} and we provide the procedure that must be followed in the case of a first order involutive system with no zero order equation, for example an involutive Spencer form.\\

\noindent 
$\bullet$  Look at the equations of class $n$ solved with respect to $y^1_n,...,y^{\beta}_n$.\\
$\bullet$  Use integrations by parts like:\\
\[ y^1_n-a(x)y^{\beta +1}_n=d_n(y^1-a(x)y^{\beta +1})+{\partial}_na(x)y^{\beta +1}={\bar{y}}^1_n+{\partial}_na(x)y^{\beta +1}  \]
$\bullet$  Modify $y^1,...,y^{\beta} $ to ${\bar{y}}^1,...,{\bar{y}}^{\beta}$ in order to "{\it absorb}" the various $y^{\beta +1}_n,...,y^m_n$ {\it only appearing in the equations of class} $n$.\\

We have the following unexpected result providing what we shall call a {\it reduced Spencer form}:\\

\noindent
{\bf THEOREM 2.4}: The new equations of class $n$ contain $y^1,...,y^{\beta}$ and their jets but only contain $y^{\beta +1}_i,...,y^m_i$ with $0\leq i\leq n-1$ while the equations of class $1,...,n-1$ no longer contain $y^{\beta+1},...,y^m$ and their jets. Accordingly, as we shall see in the next section, any torsion element, if it exists, only depends on ${\bar{y}}^1,...,{\bar{y}}^{\beta}$.\\

If ${\chi}_1, ... , {\chi}_n$ are $n$ algebraic indeterminates or, in a more intrinsic way, if $\chi={\chi}_idx^i\in T^*$ is a covector and ${\cal{D}}:E \longrightarrow F:\xi \longrightarrow a^{\tau\mu}_k(x){\partial}_{\mu}{\xi}^k(x)$ is a linear {\it involutive} operator of order $q$, we may introduce the {\it characteristic matrix} $a(x,\chi)=(a^{\tau\mu}_k(x){\chi}_{\mu}, \mid \mu\mid={\mu}_1 + ... +{\mu}_n=q)$ and the resulting map ${\sigma}_{\chi}({\cal{D}}):E \longrightarrow F$ is called the {\it symbol} of ${\cal{D}}$ at $\chi$. Then there are two possibilities:  \\
\noindent
$\bullet$ If $max_{\chi}rk({\sigma}_{\chi}({\cal{D}})< m \Leftrightarrow {\alpha}^n_q>0$: the characteristic matrix fails to be injective for any covector.\\
\noindent
$\bullet$ If $max_{\chi}rk({\sigma}_{\chi}({\cal{D}})= m\Leftrightarrow {\alpha}^n_q=0$: the characteristic matrix fails to be injective if and only if all the determinants of the $m\times m$ submatrices vanish. However, one can prove that this algebraic ideal $\mathfrak{a} \in K[\chi]$ is not intrinsically defined and must be replaced by its radical $rad(\mathfrak{a})$ made by all polynomials having a power in $\mathfrak{a}$. This radical ideal is called the {\it characteristic ideal} of the operator.\\

\noindent
{\bf DEFINITION 2.5}: For each $x\in X$, the algebraic set defined by the characteristic ideal is called the {\it characteristic set} of $\cal{D}$ at $x$ and $V={\cup}_{x\in X}V_x$ is called the {\it characteristic set} of $\cal{D}$ while we keep the word "{\it variety}" for an irreducible algebraic set defined by a prime ideal.\\

One has the following important theorem ([24], [38]) that will play an important part later on:  \\

\noindent
{\bf THEOREM 2.6}: (Hilbert-Serre) The {\it dimension} $d(V)$ of the characteristic set, that is the maximum dimension of the irreducible components, is equal to the number of non-zero characters while the {\it codimension} $cd(V)= n-d(V)$ is equal to the number of zero characters, that is to the number of "{\it full} " classes in the Janet tabular of an involutive system.   \\

If $P=a^{\mu}d_{\mu}\in D=K[d]$ with implicit summation on the multi-index, the highest value of ${\mid}\mu {\mid}$ with $a^{\mu}\neq 0$ is called the {\it order} of the {\it operator} $P$ and the ring $D$ with multiplication $(P,Q)\longrightarrow P\circ Q=PQ$ is filtred by the order $q$ of the operators. We have the {\it filtration} $0\subset K=D_0\subset D_1\subset  ... \subset D_q \subset ... \subset D_{\infty}=D$. Moreover, it is clear that $D$, as an algebra, is generated by $K=D_0$ and $T=D_1/D_0$ with $D_1=K\oplus T$ if we identify an element $\xi={\xi}^id_i\in T$ with the vector field $\xi={\xi}^i(x){\partial}_i$ of differential geometry, but with ${\xi}^i\in K$ now. It follows that $D={ }_DD_D$ is a {\it bimodule} over itself, being at the same time a left $D$-module by the composition $P \longrightarrow QP$ and a right $D$-module by the composition $P \longrightarrow PQ$. We define the {\it adjoint} map $ad:D \longrightarrow D^{op}:P=a^{\mu}d_{\mu} \longrightarrow  ad(P)=(-1)^{\mid \mu \mid}d_{\mu}a^{\mu}$ and we have $ad(ad(P))=P$. It is easy to check that $ad(PQ)=ad(Q)ad(P), \forall P,Q\in D$. Such a definition can also be extended to any matrix of operators by using the transposed matrix of adjoint operators (See [24-26],[30],[40-42] for more details and applications to control theory and mathematical physics). \\

Accordingly, if $y=(y^1, ... ,y^m)$ are differential indeterminates, then $D$ acts on $y^k$ by setting $d_{\mu}y^k=y^k_{\mu}$ with $d_iy^k_{\mu}=y^k_{\mu+1_i}$ and $y^k_0=y^k$. We may therefore use the jet coordinates in a formal way as in the previous section. Therefore, if a system of OD/PD equations is written in the form:  \\
\[ {\Phi}^{\tau}\equiv a^{\tau\mu}_ky^k_{\mu}=0\] 
with coefficients $a^{\tau \mu}_k \in K$, we may introduce the free differential module $Dy=Dy^1+ ... +Dy^m\simeq D^m$ and consider the differential submodule $I=D\Phi\subset Dy$ which is usually called the {\it module of equations}, both with the {\it differential module} $M=Dy/D\Phi$ or $D$-module and we may set $M={ }_DM$ if we want to specify the ring of differential operators. The work of Macaulay only covers the case $m=1$ with $K$ replaced by $k\subseteq cst(K)$. Again, we may introduce the formal {\it prolongation} with respect to $d_i$ by setting:  \\
\[  d_i{\Phi}^{\tau}\equiv a^{\tau\mu}_ky^k_{\mu+1_i}+({\partial}_ia^{\tau\mu}_k)y^k_{\mu}\] 
in order to induce maps $d_i:M \longrightarrow M:{\bar{y} }^k_{\mu} \longrightarrow {\bar{y}}^k_{\mu+1_i}$ if we use to denote the residue $Dy \longrightarrow M: y^k \longrightarrow {\bar{y}}^k$ by a bar as in algebraic geometry. However, for simplicity, we shall not write down the bar when the background will indicate clearly if we are in $Dy$ or in $M$.\\

As a byproduct, the differential modules we shall consider will always be {\it finitely generated} ($k=1,...,m<\infty$) and {\it finitely presented} ($\tau=1, ... ,p<\infty$). Equivalently, introducing the {\it matrix of operators} ${\cal{D}}=(a^{\tau\mu}_kd_{\mu})$ with $m$ columns and $p$ rows, we may introduce the morphism $D^p \stackrel{{\cal{D}}}{\longrightarrow} D^m:(P_{\tau}) \longrightarrow (P_{\tau}{\Phi}^{\tau}):P \longrightarrow P\Phi=P{\cal{D}}$ over $D$ by acting with $D$ {\it on the left of these row vectors} while acting with ${\cal{D}}$ {\it on the right of these row vectors} and the {\it presentation} of $M$ is defined by the exact cokernel sequence $D^p \longrightarrow D^m \longrightarrow M \longrightarrow 0 $. It is essential to notice that the presentation only depends on $K, D$ and $\Phi$ or $ \cal{D}$, that is to say never refers to the concept of (explicit or formal) solutions. It is at this moment that we have to take into account the results of the previous section in order to understant that certain presentations will be much better than others, in particular to establish a link 
with formal integrability and involution. \\

\noindent
{\bf DEFINITION 2.7}: It follows from its definition that $M$ can be endowed with a {\it quotient filtration} obtained from that of $D^m$ which is defined by the order of the jet coordinates $y_q$ in $D_qy$. We have therefore the {\it inductive limit} $0 \subseteq M_0 \subseteq M_1 \subseteq ... \subseteq M_q \subseteq ... \subseteq M_{\infty}=M$ with $d_iM_q\subseteq M_{q+1}$ and $M=DM_q$ for $q\gg 0$ with prolongations $D_rM_q\subseteq M_{q+r}, \forall q,r\geq 0$. We shall set $gr(M_q)= G_q=M_q/M_{q-1}$ and $gr(M)=G={\oplus}_qG_q$. \\

Having in mind that $K$ is a left $D$-module for the action $(D,K) \longrightarrow K:(d_i,a)\longrightarrow {\partial}_ia$ and that $D$ is a bimodule over itself, {\it we have only two possible constructions}:  \\

\noindent
{\bf DEFINITION 2.8}: We define the {\it system} $R=hom_K(M,K)=M^*$ and set $R_q=hom_K(M_q,K)=M_q^*$ as the {\it system of order} $q$. We have the {\it projective limit} $R=R_{\infty} \longrightarrow ... \longrightarrow R_q \longrightarrow ... \longrightarrow R_1 \longrightarrow R_0$. It follows that $f_q\in R_q:y^k_{\mu} \longrightarrow f^k_{\mu}\in K$ with $a^{\tau\mu}_kf^k_{\mu}=0$ defines a {\it section at order} $q$ and we may set $f_{\infty}=f\in R$ for a {\it section} of $R$. For a ground field of constants $k$, this definition has of course to do with the concept of a formal power series solution. However, for an arbitrary differential field $K$, {\it the main novelty of this new approach is that such a definition has nothing to do with the concept of a formal power series solution} ({\it care}) as illustrated in ([27]).\\

\noindent
{\bf DEFINITION 2.9}: We may define the right differential module $hom_D(M,D)$.  \\

\noindent
{\bf PROPOSITION 2.10}: When $M$ is a left $D$-module, then $R$ is also a left $D$-module. \\

\noindent
{\it Proof}: As $D$ is generated by $K$ and $T$ as we already said, let us define:  \\
\[  (af)(m)=af(m), \hspace{4mm} \forall a\in K, \forall m\in M \]
\[ (\xi f)(m)=\xi f(m)-f(\xi m), \hspace{4mm} \forall \xi=a^id_i\in T,\forall m\in M  \]
In the operator sense, it is easy to check that $d_ia=ad_i+{\partial}_ia$ and that $\xi\eta - \eta\xi=[\xi,\eta]$ is the standard bracket of vector fields. We finally 
get $(d_if)^k_{\mu}=(d_if)(y^k_{\mu})={\partial}_if^k_{\mu}-f^k_{\mu +1_i}$ and thus recover {\it exactly} the Spencer operator though {\it this is not evident at all}. We also get $(d_id_jf)^k_{\mu}={\partial}_{ij}f^k_{\mu}-{\partial}_if^k_{\mu+1_j}-{\partial}_jf^k_{\mu+1_i}+f^k_{\mu+1_i+1_j} \Longrightarrow d_id_j=d_jd_i, \forall i,j=1,...,n$ and thus $d_iR_{q+1}\subseteq R_q\Longrightarrow d_iR\subset R$ induces a well defined operator $R\longrightarrow T^*\otimes R:f \longrightarrow dx^i\otimes d_if$. This result has been discovered (up to sign) by Macaulay  in 1916 ([15]). For more details on the Spencer operator and its applications, the reader may look at ([21],[22],[28],[48]).  \\
\hspace*{12cm}   Q.E.D.  \\

\noindent
{\bf DEFINITION 2.11}: $t_r(M)$ is the greatest differential submodule of $M$ having codimension $> r$.  \\

\noindent
{\bf PROPOSITION 2.12}: $cd(M)=cd(V)=r \Longleftrightarrow {\alpha}^{n-r}_q\neq 0, {\alpha}^{n-r+1}_q= ... ={\alpha}^n_q=0 \Longleftrightarrow t_r(M)\neq M, t_{r-1}(M)= ... =t_0(M)=t(M)=M$ and this intrinsic result can be most easily checked by using the standard or reduced Spencer form of the system defining $M$.  \\

We are now in a good position for defining and studying purity for differential modules. \\

\noindent
{\bf DEFINITION 2.13}: $M$ is $r$-{\it pure} $\Longleftrightarrow t_r(M)=0, t_{r-1}(M)=M \Longleftrightarrow cd(Dm)=r, \forall m\in M$. More generally, $M$ is {\it pure} if it is $r$-pure for a certain $0\leq r\leq n$ and $M$ is {\it pure} if it is $r$-pure for a certain $0\leq r \leq n$. In particular, $M$ is $0$-pure if $t(M)=0$ and, if $cd(M)=r$ but $M$ is not $r$-pure, we may call $M/t_r(M)$ the {\it pure part} of $M$. It follows that $t_{r-1}(M)/t_r(M)$ is equal to zero or is 
$r$-pure (See the picture in [20], p 545). When $M=t_{n-1}(M)$ is $n$-pure, its defining system is a finite dimensional vector space over $K$ with a symbol of finite type, that is when $g_q=0$ is (trivially) involutive. Finally, when $t_{r-1}(M)=t_r(M)$, we shall say that there is a {\it gap} in the purity filtration:   \\
\[   0=t_n(M) \subseteq t_{n-1}(M) \subseteq ... \subseteq t_1(M) \subseteq t_0(M)=t(M) \subseteq M     \]

\noindent
{\bf PROPOSITION 2.14}: $t_r(M)$ does not depend on the presentation or on the filtration of $M$.  \\

\noindent
{\bf EXAMPLE 2.15}: If $K=\mathbb{Q}$ and $M$ is defined by the involutive system $y_{33}=0, y_{23}=0, y_{13}=0$, then $z=y_3$ satifies $d_3z=0, d_2z=0, d_1z=0$ and $cd(Dz)=3$ while $z'=y_2$ only satisfies $d_3z'=0$ and $cd(Dz')=1$. We have the purity filtration 
$0 =t_3(M) \subset t_2(M) =t_1(M) \subset t_0(M)=t(M)=M$ with one gap and two strict inclusions.  \\

We now recall the definition of the {\it extension modules} $ext_D^i(M,D)$ that we shall simply denote by $ext^i(M)$ and the way to use their dimension or codimension. We point out once more that these numbers can be most easily obtained by bringing the underlying systems to involution in order to get informations on $M$ from informations on $G$. We divide the procedure into four steps that can be achieved by means of computer algebra ([43],[44]): \\

\noindent
$\bullet$ Construct a {\it free resolution} of $M$, say:  \\
\[   ... \longrightarrow F_i \longrightarrow ... \longrightarrow F_1 \longrightarrow F_0 \longrightarrow M \longrightarrow 0  \]

\noindent
$\bullet$ Suppress $M$ in order to obtain the {\it deleted sequence}:  \\
\[     ... \longrightarrow F_i \longrightarrow ... \longrightarrow F_1 \longrightarrow F_0 \longrightarrow 0  \hspace{11mm}  \]

\noindent
$\bullet$ Apply $hom_D(\bullet,D)$ in order to obtain the {\it dual sequence} heading backwards: \\
\[     ... \longleftarrow hom_D(F_i,D) \longleftarrow ... \longleftarrow hom_D(F_1,D) \longleftarrow hom_D(F_0,D) \longleftarrow 0    \]

\noindent
$\bullet$ Define $ext^i(M)$ to be the cohomology at $hom_D(F_i,D)$ in the dual sequence in such a way that $ext^0(M)=hom_D(M,D)$.  \\

The following nested chain of difficult propositions and theorems can be obtained, {\it even in the non-commutative case}, by combining the use of extension modules and {\it bidualizing complexes} in the framework of algebraic analysis. The main difficulty is to obtain first these results for the {\it graded module} $G=gr(M)$ by using techniques from commutative algebra before extending them to the {\it filtred module} $M$ as in ([1],[2],[10],[13],[14],[16],[24],[29],[40],[48]).  \\

\noindent
{\bf THEOREM 2.16}: The extension modules do not depend on the resolution of $M$ used.  \\

\noindent
{\bf PROPOSITION 2.17}: Applying $hom_D(\bullet,D)$ provides right $D$-modules that can be transformed to left $D$-modules by means of the {\it side changing functor} and vice-versa. Namely, if $N_D$ is a right $D$-module, then ${}_DN={\wedge}^nT{\otimes}_KN_D$ is the {\it converted left} $D$-module while, if $N={}_DN$ is a left $D$-module, then $N_D={\wedge}^nT^*{\otimes}_KN$ is the {\it converted right} $D$-module.\\

\noindent
{\bf PROPOSITION 2.18}: Instead of using $hom_D(\bullet,D)$ and the side changing functor in the module framework, we may use $ad$ in the operator framework. Namely, to any operator ${\cal{D}}:E \longrightarrow F$ we may associate the formal adjoint $ad({\cal{D}}):{\wedge}^nT^*\otimes F^*\longrightarrow {\wedge}^nT^*\otimes E^*$ with the useful though striking relation $rk_D(ad({\cal{D}}))=rk_D({\cal{D}})$.  \\

\noindent
{\bf PROPOSITION 2.19}: $ext^i(M)$ is a torsion module $\forall 1\leq i \leq n$ but $ext^0(M)=hom_D(M,D)$ may not be a torsion module.  \\

\noindent
{\bf EXAMPLE 2.20}: When $M$ is a torsion module, we have $hom_D(M,D)=0$ (exercise). When $n=3$ and the torsion-free module $M$ is defined by the 
formally surjective $div$ operator, the formal adjoint of $div$ is $-grad$ which defines a torsion module. Also, when $n=1$ as in classical control theory, a controllable system with coefficients in a differential field allows to define a torsion-free module $M$ which is free in that case because a finitely generated module over a principal ideal domain is free if and only if it is torsion-free and $hom_D(M,D)$ is thus also a free module. \\

\noindent
{\bf THEOREM 2.21}: \hspace{1cm}  $ext^i(M)=0, \forall i<cd(M)$ and $\forall i\geq n+1$.  \\

\noindent
{\bf THEOREM 2.22}: \hspace{1cm}  $cd(ext^i(M))\geq i$.  \\

\noindent
{\bf THEOREM 2.23}: \hspace{1cm}  $cd(M)\geq r \Leftrightarrow ext^i(M)=0, \forall i<r$.  \\

\noindent
{\bf PROPOSITION 2.24}: \hspace{1cm}  $cd(M)=r \Longrightarrow cd(ext^r(M))=r$ and $ext^r(M)$is $r$-pure.  \\

\noindent
{\bf PROPOSITION 2.25}: \hspace{1cm}  $ext^r(ext^r(M))$ is equal to $0$ or is $r$-pure, $\forall 0\leq r \leq n$.  \\

\noindent
{\bf PROPOSITION 2.26}: If we set $t_{-1}(M)=M$, there are exact sequences $\forall 0\leq r \leq n$:   \\
\[       0 \longrightarrow t_r(M) \longrightarrow t_{r-1}(M)  \longrightarrow ext^r(ext^r(M))  \]

\noindent
{\bf THEOREM 2.27}: If $cd(M)=r$, then $M$ is $r$-pure if and only if there is a monomorphism $0 \longrightarrow M \longrightarrow ext^r(ext^r(M))$ of left differential modules.  \\

\noindent
{\bf THEOREM 3.28}: $M$ is pure $\Longleftrightarrow ext^s(ext^s(M))=0 , \forall s\neq cd(M)$.  \\

\noindent
{\bf COROLLARY 2.29}: If $M$ is $r$-pure with $r\geq 1$, then it can be embedded into a differential module $L$ having a free resolution with only $r$ operators.   \\

The previous theorems are known to characterize purity but it is however evident that they are not very useful in actual practice. For more details on these two results which are absolutely out of the scope of this paper, see ([2], p 490-491) and ([24], p 547). Proposition 3.24 and Theorem 3.25 come from the  Cohen-Macaulay property of $M$, namely $cd(M)=g(M)= inf \{ i\mid ext^i(M)\neq 0\}$ where $g(M)$ is called the {\it grade} of $M$ (See [2]  and [24],[29] for more details). \\

\noindent
{\bf THEOREM 2.30}: When $M$ is $r$-pure, the characteristic ideal is thus {\it unmixed}, that is a finite intersection of prime ideals having the same codimension $r$ and the characteristic set is {\it equidimensional}, that is the union of irreducible algebraic varieties having the same codimension $r$.  \\

In $2013$ we have provided a new effective test for checking purity while using the involutivity of the Spencer form with four steps as follows ([29]):  \\

\noindent
$\bullet$  STEP 1: Compute the involutive Spencer form of the system and the number $r$ of full classes.  \\  

\noindent
$\bullet$  STEP 2: Select only the equations of class $1$ to $d(M)=n-r$ of this Spencer form which are making an involutive system over $K[d_1,...,d_{(n-r)}]$.  \\  

\noindent
$\bullet$  STEP 3: Using differential biduality for such a system, check if it defines a torsion-free module $M_{(n-r)}$ and work out a parametrization.   \\ 
 
\noindent
$\bullet $ STEP 4: Substitute the above parametrization in the remaning equations of class $n-r+1, ... , n$ of the Spencer form in order to get a system of PD equations which provides the {\it parametrizing module} $L$ in such a way that $M\subseteq L$ and $L$ has a resolution with $r$ operators.  \\

\noindent
{\bf THEOREM 2.31}: As purity is an intrinsic property, we may work with an involutive Spencer form and $M$ is $r$-pure if the classes $n-r+1,...,n$ are full and the module $M_{(n-r)}$ defined by the equations of class $1+...+$ class $(n-r)$ is torsion-free. Hence $M$ is $0$-pure if it is torsion-free. \\

We shall now illustrate and apply this new procedure in the next two sections.  \\

\noindent
{\bf 3) MOTIVATING EXAMPLES}:  \\

\noindent
{\bf EXAMPLE 3.1}: With $n=3, m=1$ and $K=\mathbb{Q}$, let us consider the following polynomial ideal:  \\
\[   \mathfrak{a}= (({\chi}_3)^2, {\chi}_2{\chi}_3 - {\chi}_1{\chi}_3, ({\chi}_2)^2 - {\chi}_1{\chi}_2) \subset K[{\chi}_1,{\chi}_2,{\chi}_3]=K[\chi] \]
We shall discover that it is not evident to prove that it is an unmixed polynomial ideal and that the corresponding differential module is $1$-pure.  \\
The first result is provided by the existence of the {\it primary decomposition} obtained from the two existing factorizations:  \\
\[    \mathfrak{a}= (({\chi}_3)^2, {\chi}_2 - {\chi}_1) \cap ({\chi}_3,{\chi}_2)={\mathfrak{q}}' \cap {\mathfrak{q}}"  \]
Taking the respective radical ideals, we get the {\it prime decomposition}:  \\
\[  rad(\mathfrak{a})= ({\chi}_3, {\chi}_2 - {\chi}_1) \cap ({\chi}_3,{\chi}_2)= {\mathfrak{p}}' \cap {\mathfrak{p}}" = 
rad({\mathfrak{q}}') \cap rad({\mathfrak{q}}")  \]
The corresponding involutive system is:   \\
\[  \left\{  \begin{array}{lcl}
 y_{33}  & = &0  \\
 y_{23} - y_{13} & = & 0   \\
 y_{22}-y_{12} & = & 0
\end{array}  
\right.  \fbox{$\begin{array}{lll}
1 & 2 & 3     \\
1 & 2 & \bullet  \\
1 & 2 & \bullet  
\end{array}  $  }    \]
with characters $({\alpha}^3_2=1-1=0, {\alpha}^2_2=2-2=0, {\alpha}^1_2=3-0=3$ and $dim(g_2)=\sum \alpha=3$.  \\
Setting $(z^1=y,z^2=y_1,z^3=y_2,z^4=y_3)$, we obtain the involutive first order Spencer form:  \\
\[  \left\{  \begin{array}{l}
   z^4_3  =0, z^3_3-z^4_1=0,z^2_3- z^4_1=0, z^1_3 - z^4= 0 \\
  z^4_2 - z^4_1=0, z^3_2 - z^3_1=0, z^2_2- z^3_1=0, z^1_2 - z^3=0  \\ 
 z^1_1 - z^2  =  0
\end{array}  
\right.  \fbox{$\begin{array}{lll}
1 & 2 & 3     \\
1 & 2 & \bullet  \\
1 & \bullet & \bullet  
\end{array}  $  }    \]
with new characters ${\alpha}^3_1=4-4=0, {\alpha}^3_1=4-4=0, {\alpha}^1_1=4-1=3$ and similarly $dim(g_1)=\sum \alpha=3$.
Both class $3$ and class $2$ are full while class $1$ is defining a torsion-free module $M_{(1)}$ over $K[d_1]$ by means of a trivially involutive system of class $1$. Hence the differential module $M$ is such that $cd(M)=2$ and is $1$-pure because it is $1$-pure in {\it this} presentation.  \\
Suppressing the bar for the various residues, we are ready to exhibit the {\it relative parametrization} defining the {\it parametrization module} $L$ because we may choose the $3$ potentials $(z^1=y, z^3, z^4)$ while taking into account that $z^2=y_1=d_1y$:  \\
\[  \left\{   \begin{array}{lcl}
 z^4_3   & = &   0 \\
    z^3_3-z^4_1& = & 0 \\
      z^1_3 - z^4& = & 0 \\
  z^4_2 - z^4_1& = & 0 \\
  z^3_2 - z^3_1 & = & 0  \\
 z^1_2 - z^3 & = & 0 
\end{array} \right.  \fbox{$ \begin{array}{ccc}
1 & 2 & 3 \\
1 & 2 & 3  \\
1 & 2 & 3   \\
1 & 2 & \bullet  \\
1 & 2 & \bullet  \\
1 & 2 & \bullet
\end{array} $ }  \]
Both $(y, z^3,z^4)$ are torsion elements and we can eliminate $(z^3,z^4)$ in order to find the desired system that must be satisfied by $y$ which is showing the inclusion $M\subset L$ but we have indeed $M=L$ because $z^3=y_2, z^4=y_3$. It follows that $M$ admits a free resolution with only $2$ operators, a result following at once from the last Janet tabular, {\it contrary to the previous one}. \\
The reader may treat similarly the example $\mathfrak{a}= ({\chi}_1,{\chi}_2) \cap ({\chi}_3, {\chi}_4)$ and look at ([27]) for details. (Hint: use the involutive system $y_{44}+y_{14}=0, y_{34} +y_{13}=0,y_{33}+y_{23}=0, y_{24}-y_{13}=0$).   \\

\noindent
{\bf EXAMPLE 3.2}: With $n=3, m=1, q=2, K=\mathbb{Q}, D=K[d_1,d_2,d_3]$, let us consider the differential module $M$ defined by the second order system $Py\equiv y_{33}=0, Qy\equiv y_{13} - y_2=0$ first considered by Macaulay in $1916$ ([15],[39]). We shall prove that $M$ is $2$-pure through the inclusion $0\rightarrow M \rightarrow ext^2(ext^2(M))$ directly and by finding out a relative parametrization, a result highly not evident at first sight. \\
First of all, in order to find out the codimension $cd(M)=2$, we have to consider the equivalent involutive system:  \\
\[  \left\{  \begin{array}{lclcl}
{\Phi}^4 & \equiv & y_{33}  & = &u  \\
{\Phi}^3 & \equiv & y_{23}  & = & d_1u-d_3v   \\
{\Phi}^2 & \equiv & y_{22} & = & d_{11}u-d_{13}v-d_2v \\
{\Phi}^1 & \equiv & y_{13}-y_2 & = & v
\end{array}  
\right.  \fbox{$\begin{array}{lll}
1 & 2 & 3     \\
1 & 2 & \bullet  \\
1 & 2 & \bullet  \\
1 & \bullet & \bullet 
\end{array}  $  }    \]
The Janet tabular on the rigt allows at once to compute the characters ${\alpha}^3_2=0,{\alpha}^2_2=0, {\alpha}^1_2=\alpha=3-1=2$ and to construct the following strictly exact sequence of differential modules:  \\
\[   0 \rightarrow D \rightarrow D^4 \rightarrow D^4 \rightarrow D \stackrel{p}{ \rightarrow} M \rightarrow 0  \]
Also, we have $rad(\mathfrak{a})=rad (({\chi}_3)^2, {\chi}_2{\chi}_3, ({\chi}_2)^2, {\chi}_1{\chi}_3) = ({\chi}_3, {\chi}_2)= \mathfrak{p}\Rightarrow dim(V)=1$.  \\
As the classes $3$ and $2$ are full, it follows that $d(M)=d(Dy)=1\Rightarrow cd(M)=n-1=2$ if we denote simply by $y$ the canonical residue $\bar{y}$ of $y$ after identifyoing $D$ with $Dy$. We have constructed explicitly in ([29]) a finite length resolution of $N=ext^2(M)$ by pointing out that $N$ does not depend on the resolution of $M $ used and one can refer to the single compatibility condition (CC) $P\circ Q y - Q\circ P y=0$ for the initial system in the exact sequence made by second order operators:  \\
\[  \hspace{1cm} 0 \longrightarrow D  \underset {2}{\stackrel{{\cal{D}}_1}{\longrightarrow}} D^2  \underset {2}{\stackrel{\cal{D}}{\longrightarrow}} D  \stackrel{p}{\longrightarrow} M \longrightarrow 0  \]
Indeed, introducing differential duality through the functor $hom_D(\bullet,D)$ and the respective adjoint operators, we may define the torsion left differential module $N$ by the long exact sequence:  \\
\[  0 \longleftarrow N \stackrel{q}{\longleftarrow} D  \stackrel{ad({\cal{D}}_1)}{\longleftarrow} D^2 \stackrel{ad(\cal{D})}{\longrightarrow} D   \longleftarrow 0 \hspace{12mm}  \]
showing that $rk_D(M)=0 \Rightarrow rk_D(N)=1-2+1=0$ because of the additivity property of the differential rank and the vanishing of the Euler-Poincar\'{e} characteristic of the full sequence. It follows that $M={ext}^2(N)=ext^2(ext^2(M))$.   \\
Similarly, using certain parametric jet variables as new unknowns, we may set $z^1=y, z^2=y_1, z^3=y_2, z^4=y_3$ in order to obtain the following involutive first order system with no zero order equation:  \\
\[  \left \{ \begin{array}{lcc}
class \hspace{1mm}  3    & \hspace{1cm}                 &       d_3z^1-z^4=0, d_3z^2-z^3=0, d_3z^3=0,d_3z^4=0  \\
class \hspace{1mm}  2    & \hspace{1cm}                 &       d_2z^1-z^3=0, d_2z^2-d_1z^3=0, d_2z^3=0, d_2z^4=0  \\
class \hspace{1mm}  1    & \hspace{1cm}                 &                      d_1z^1-z^2=0, d_1z^4-z^3=0   
 \end{array} 
\right.   \fbox{  $   \begin{array}{lll}
1  & 2  &  3  \\
1  &  2 &  \bullet  \\
1  &  \bullet  &  \bullet  
\end{array}  $  }   \]

\noindent
where we have separated the classes while using standard computer algebra notations this time instead of the jet notations used in the previous example. Contrary to what could be believed, this operator does not describe the Spencer sequence that could be obtained from the previous Janet sequence but we can use it exactly like a Janet sequence or exactly like a Spencer sequence. We obtain therefore a long strictly exact sequence of differential modules {\it with only first order operators} while replacing $Dy$ by $Dz=Dz^1 +Dz^2 + Dz^3 +Dz^4$ as follows:  \\
\[   0 \rightarrow D ^2  \underset 1{\rightarrow} D^8 \underset 1{\rightarrow} D^{10} \underset 1{\rightarrow} D^4 \rightarrow M \rightarrow 0  \] 
and we still have the vanishing Euler-Poincar\'{e} characteristic $ 2 - 8 +10 - 4 =0$.  \\
The differential module $M_1$ is defined over $K[d_1]$ by the two PD equations of class $1$ and is easily seen to be torsion-free with the two potentials $(z^1=y, z^4)$. Substituting into the PD equations of class $2$ and $3$, we obtain the generating differential constraints: \\
\[   \left\{  \begin{array}{l} 
d_3z^1 - z^4=0  \\
 d_3z^4=0 \\
 d_2z^1 - d_1z^4=0  \\
  d_2z^4=0  
\end{array}  \right.  \fbox{ $ \begin{array}{ccc}
1 & 2 & 3  \\
1 & 2 & 3   \\  
 1 & 2 & \bullet  \\
 1 & 2 & \bullet  
 \end{array} $  }  \]
They define the {\it parametrization module} $L$ and the inclusion $M\subseteq L$ is obtained by eliminating $z^4$ but we have indeed $M=L$ because 
$z^4=d_3y$.  \\

\noindent
{\bf EXAMPLE 3.3}: We have provided in ([29], Example 4.2) a case leading to a strict inclusion $M \subset L$ that we revisit now totally in this new framework. With $K=\mathbb{Q}, m=1, n=4, q=2$, let us study the $2$-pure differential module $M$ defined by the involutive second order system:   \\
\[  \left\{  \begin{array}{lcl}
 y_{44}  & = & 0  \\
 y_{34}  & = &  0   \\
 y_{33} & = &  0  \\
 y_{24}-y_{13} & =&   0
\end{array}  
\right.  \fbox{$\begin{array}{llll}
1 & 2 & 3  & 4   \\
1 & 2 & 3 &\bullet  \\
1 & 2 & 3 &\bullet  \\
1 &2 &\bullet & \bullet 
\end{array}  $  }    \]

\noindent
From the Janet tabular we may construct at once the Janet sequence:  \\
\[  0 \longrightarrow \Theta \longrightarrow {\:1\:} \stackrel{\cal{D}}{\longrightarrow} {\:4\:} \stackrel{{\cal{D}}_1}{\longrightarrow} {\:4\:}
\stackrel{{\cal{D}}_2}{\longrightarrow} {\:1\:} \longrightarrow 0   \]

\noindent
where ${\cal{D}}_1$ is defined by the involutive system:  \\
\[  \left\{  \begin{array}{lcl}
 d_4(y_{34})-d_3(y_{44})  & =    &0  \\
 d_4(y_{33})-d_3(y_{34})& =  & 0   \\
 d_4(y_{24}-y_{13})-d_2(y_{44})+d_1(y_{34}) & =  & 0  \\
 d_3(y_{24} - y_{13})-d_2(y_{34})+d_1(y_{33}) & = & 0
\end{array}  
\right.  \fbox{$\begin{array}{llll}
1 & 2 & 3  & 4   \\
1 & 2 & 3 & 4 \\
1 & 2 & 3 & 4 \\
1 & 2 & 3& \bullet 
\end{array}  $  }    \]
and so on. We have therefore a free resolution of $M$ with $3$ operators:  \\
\[ 0 \longrightarrow D \longrightarrow D^4 \longrightarrow D^4 \longrightarrow D \longrightarrow M  \longrightarrow 0  \]
and thus discover that $pd(M)\leq 3$. \\
However, we have $rad(({\chi}_4)^2, {\chi}_3{\chi}_4, ({\chi}_3)^2, {\chi}_2{\chi}_4 - {\chi}_1{\chi}_3)=({\chi}_4,{\chi}_3)= \mathfrak{p} \Rightarrow dim (V)=1$.   \\
Let us transform the initial second order involutive system for $y$ into a first order involutive system for $(z^1=y, z^2=y_1,z^3=y_2, z^4=y_3,z^5=y_4)$ as follows:  \\
\[ \left\{ \begin{array}{lcl}
d_4z^1-z^5=0, d_4z^2-d_1z^5=0, d_4z^3-d_1z^4=0, d_4z^4=0, d_4z^5=0 & &   \\
d_3z^1-z^4=0, d_3z^2-d_1z^4=0, d_3z^3-d_2z^4=0, d_3z^4=0, d_3z^5=0 & &  \\
d_2z^1-z^3=0, d_2z^2-d_1z^3=0, d_2z^5-d_1z^4=0 &  &   \\
d_1z^1-z^2=0 &  & 
\end{array}
\right.  \fbox{ $\begin{array}{llll}
1 & 2 & 3 & 4 \\
1 & 2 & 3 & \bullet \\
1 & 2 & \bullet & \bullet \\
1 & \bullet & \bullet & \bullet 
\end{array} $ }  \]
with five equations of full class $4$,five equations of full class $3$, three equations of class $2$ and finally one equation of class $1$. The equations of classes 2 and 1 are providing an involutive system over $\mathbb{Q}[d_1,d_2]$ defining a torsion-free module $M_{(2)}$ that can be parametrized by setting $z^1=y, z^2=d_1y, z^3=d_2y, z^4=d_2z,z^5=d_1z$ with only $2$ arbitrary potentials $(y,z)$. Substituting in the other equations of classes 3 and 4, we finally discover that $L$ is defined by the involutive system describing the relative parametrization:  \\
\[  \left\{ \begin{array}{lcl}
 d_4y-d_1z & =  & 0 \\
d_4z      & = & 0 \\
d_3y-d_2z & = & 0 \\
 d_3z   & = & 0  
\end{array}
\right.  \fbox{ $ \begin{array}{llll}
1 & 2 & 3 & 4  \\
1 & 2 & 3 & 4 \\
1 & 2 & 3 & \bullet \\
1 & 2 & 3 & \bullet 
\end{array} $} \]
We have the strict inclusion $M \subset L$ obtained by eliminating $z$ because now $z \notin Dy$ if we take the residue or, equivalently, the residue of $z$ does not belong to $M$. The differential module $L$ defined by the above system is therefore $2$-pure with a strict inclusion $M\subset L$ and admits a free resolution with only $2$ operators according to its Janet tabular. \\

\noindent
{\bf EXAMPLE 3.4}: ({\it Contact structure}) With $n=m=3$ and $K=\mathbb{Q}(x^1,x^2,x^3)$ let us introduce the so-called {\it contact} $1$-form $\alpha=dx^1 - x^3dx^2$ and consider the first order system of infinitesimal Lie equations obtained by eliminating the contact factor $\rho$ from the equations ${\cal{L}}(\xi)\alpha = \rho \alpha$. We let the reader check that he will obtain only the two equations ${\Phi}^1=0,{\Phi}^2=0$ which is nevertheless neither formally integrable nor even involutive. Using crossed derivatives one obtains the involutive system:  \\
\[  \left\{ \begin{array}{lclc}
{\Phi}^3 &\equiv & {\partial}_3{\xi}^3 + {\partial}_2 {\xi}^2 + 2 x^3 {\partial}_1{\xi}^2 & = 0  \\
{\Phi}^2 & \equiv & {\partial}_3{\xi}^1 - x^3 {\partial}_3{\xi}^2  & =   0  \\
{\Phi}^1 &  \equiv & {\partial}_2{\xi}^1 - x^3 {\partial}_2{\xi}^2 +x^3 {\partial}_1{\xi}^1 - (x^3)^2{\partial}_1{\xi}^2 - {\xi}^3 & =  0
\end{array}  \right.  \fbox{$ \begin{array}{lll}
1 & 2 & 3  \\
1 & 2 & 3  \\
1 & 2 & \bullet
\end{array} $ }  \]
with the unique CC $\Psi \equiv {\partial}_3 {\Phi}^1 - {\partial}_2{\Phi}^2 - x^3 {\partial}_1{\Phi}^2 + {\Phi}^3 =0$.
The following injective absolute parametrization is well known and we let the reader find it by using differential double duality:  \\
\[  \phi - x^3 {\partial}_3 \phi  = {\xi}^1, \,\,\, -{\partial}_3 \phi = {\xi}^2, \,\,\, {\partial}_2\phi + x^3 {\partial}_1\phi={\xi}^3 \Rightarrow 
{\xi}^1 - x^3 {\xi}^2=\phi \]
We obtain the Janet sequence
\[   \begin{array}{rcccccccl}
0 \rightarrow  & 1 & \stackrel{{\cal{D}}_{-1}}{\longrightarrow} & 3 & \stackrel{{\cal{D}}}{\longrightarrow } & 3 & \stackrel{{\cal{D}}_1}{\longrightarrow} & 1 & \rightarrow 0   \\
       &  \phi & & \xi & & \Phi & & \Psi & 
       \end{array}   \]
with formally exact adjoint sequence:  \\
 \[   \begin{array}{rcccccccl}
0 \leftarrow  & 1 & \stackrel{ad({\cal{D}}_{-1})}{\longleftarrow} & 3 & \stackrel{ad({\cal{D}})}{\longleftarrow } & 3 & \stackrel{ad({\cal{D}}_1)}{\longleftarrow} & 1 & \rightarrow 0   \\
   & \theta & & \nu & & \mu & & \lambda &   
       \end{array}   \]      
and the resolution of the trivially torsion-free module $M\simeq D$:  \\
\[  0 \longrightarrow D \longrightarrow D^3 \longrightarrow D^3 \longrightarrow M  \longrightarrow 0  \]
which splits totally because it is made with free and thus projective modules.   \\

\noindent
{\bf EXAMPLE 3.5}: ({\it Unimodular contact structure})  With $n=m=3$ and $K=\mathbb{Q}(x^1,x^2,x^3)$ let us introduce the $1$-form 
$\omega=dx^1 - x^3dx^2$ used as a geometric object and consider the first order system of infinitesimal Lie equations from the equations ${\cal{L}}(\xi)\omega = 0$. One obtains the system using jet notations:  \\
\[  {\xi}^1_3 -x^3{\xi}^2_3 =0, \,\,\, {\xi}^1_2 - x^3 {\xi}^2_2 - {\xi}^3=0, \,\, \, {\xi}^1_1 - x^3 {\xi}^2_1=0  \]
We let the reader prove that these three PD equations are differentially independent and we obtain the free resolution of $M$:  \\
\[ \hspace{1cm}  0 \longrightarrow D^3 \stackrel{{\cal{D}}}{\longrightarrow} D^3 \longrightarrow M \longrightarrow 0  \]
and its adjoint sequence is:  \\
\[  0 \longleftarrow N \longleftarrow  D^3 \stackrel{ad({\cal{D}})}{\longleftarrow} D^3 \longleftarrow 0  \hspace{13mm}  \]
because $rk_D(M=rk_D(N)= 3-3=0$, that is both $M$ and $N$ are torsion modules with $N=ext^1(M) \Rightarrow M=ext^1(N)=ext^1(ext^1(M))$ and $M$ is surely $1$-pure.
However, this system is not formally integrable, as it can be checked directly through crossed derivatives or by noticing that ${\cal{L}}(\xi)d\omega=0$ with $d\omega= dx^2 \wedge dx^3$ and ${\cal{L}}(\xi)(\omega \wedge d\omega)=0$ with $\omega \wedge d\omega = dx^1 \wedge dx^2 \wedge dx^3$. Hence, we have to add the 3 first order equations:  \\
\[        {\xi}^2_2 + {\xi}^3_3, \,\,\, {\xi}^3_1=0, \,\,\, {\xi}^2_1=0 \Rightarrow {\xi}^1_1=0 \]
Exchanging $x^1$ and $x^3$, we obtain the equivalent involutive system in $\delta$-regular coordinates:  \\
\[   \left\{  \begin{array}{lcc}
{\xi}^3_3 & = & 0  \\
{\xi}^2_3 & = & 0  \\
{\xi}^1_3   & = & 0   \\
{\xi}^2_2 + {\xi}^3_1 & = & 0 \\
{\xi}^1_2 + x^1 {\xi}^3_1 - {\xi}^3 & = & 0 \\
{\xi}^1_1 - x^1{\xi}^2_1 & = & 0
\end{array}  \right. \fbox{ $ \begin{array}{ccc}
1 & 2 & 3  \\
1 &  2 & 3  \\
1 & 2 & 3 \\
1 & 2 & \bullet  \\
1 & 2 & \bullet  \\
1 & \bullet & \bullet
\end{array} $ }  \]
The differential module $M_{(2)}$ over $K[d_1,d_2]$ is defined by the three bottom equations. Setting now $\phi= {\xi}^1 - x^1 {\xi}^2$, we deduce from the last bottom equation that ${\xi}^2=-d_1\phi$ and thus ${\xi}^1=\phi - x^1 d_1 \phi$. Finally, substituting in the equation before the last, we get ${\xi}^3 = d_2 \phi$. We have thus obtained an injective parametrization of $M_{(2)}$ which is therefore torsion-free and $M$ is $2$-pure in a coherent way. Substituting into the three upper equations, we obtain the desired relative parametrization by adding the differential constraint $d_3\phi=0$. Coming back to the original coordinates, we obtain the relative parametrization:  \\
\[      \phi - x^3 d_3\phi = {\xi}^1, \,\,\,  - d_3 \phi = {\xi}^2, \,\,\, d_2\phi= {\xi}^3  \,\,\,\, with  \,\,\,\, d_1\phi=0  \]
which is thus strikingly obtained from the previous contact parametrization by adding the only differential constraint $d_1\phi=0$.  \\   \\

\noindent
{\bf 4) APPLICATIONS}  \\  
Before studying applications to mathematical physics, we shall start with an example describing in an explicit way the Janet and Spencer sequences used thereafter, both with their link, namely the relations existing between the dimensions of the respective Janet and Spencer bundles.  \\

\noindent
{\bf EXAMPLE 4.1}: When $n=m=2,q=2$, $\omega$ is the Euclidean metric of $X={\mathbb{R}}^2$ with Christoffel symbols $\gamma$ and metric density $\tilde{\omega}=\omega / \sqrt{det(\omega)}$, we consider the two involutive systems of linear infinitesimal Lie equations $R_2 \subset {\tilde{R}}_2 \subset J_2(T)$ respectively defined by $\{{\cal{L}}(\xi)\omega =0, {\cal{L}}(\xi)\gamma=0\}$ and $\{  {\cal{L}}(\xi)\tilde{\omega}=0, {\cal{L}}(\xi)\gamma=0 \} $. We have $g_2={\tilde{g}}_2=0$ and construct the following successive commutative and exact diagrams followed by the corresponding dimensional diagrams that are used in order to construct effectively the respective Janet and Spencer differential sequences while comparing them.   \\

\[  \begin{array}{rcccccccl}
& & 0 & & 0 & & 0 & &   \\
& & & \searrow &\downarrow & &\downarrow & &   \\
& & 0 &\rightarrow &S_2T^*\otimes T&\rightarrow & F"_0 &\rightarrow  &  0 \\
& & \downarrow & & \downarrow  & \searrow & \downarrow & &  \\
0 & \rightarrow & R_2  & \rightarrow &J_2(T) &\rightarrow  & F_0& \rightarrow & 0 \\
& & \downarrow& &\downarrow & &\downarrow & &   \\
0 & \rightarrow & R_1 & \rightarrow & J_1(T)& \rightarrow & {F'}_0 & \rightarrow & 0  \\
& & \downarrow & &\downarrow  & & \downarrow &   \\
& & 0 & & 0 & &  0  & &    
\end{array}     \]

\[  \begin{array}{rcccccccl}
& & 0 & & 0 & & 0 & &   \\
& & & \searrow &\downarrow & &\downarrow & &   \\
& & 0 &\rightarrow & 6 &\rightarrow & 6 &\rightarrow  &  0 \\
& & \downarrow & & \downarrow  & \searrow & \downarrow & &  \\
0 & \rightarrow & 3  & \rightarrow & 12 &\rightarrow  & 9 & \rightarrow & 0 \\
& & \downarrow& &\downarrow & &\downarrow & &   \\
0 & \rightarrow & 3 & \rightarrow & 6 & \rightarrow & 3 & \rightarrow & 0  \\
& & \downarrow & &\downarrow  & & \downarrow &   \\
& & 0 & & 0 & &  0  & &    
\end{array}     \]

In the present situation we notice that $R_2={\rho}_1(R_1)=J_1(R_1)\cap J_2(T)\subset J_1(J_1(T))$ and thus $F_0\simeq J_1({F'}_0)$ with 
$F"_0\simeq T^*\otimes {F'}_0\simeq S_2T^*\otimes T$ by counting the dimensions because we have surely $F_0 \subset J_1({F}'_0)$ with $g_2=0$. \\

\[  \begin{array}{rcccccccccl}
& & & & 0 & & 0 & & 0  & &  \\
& & & &\downarrow & &\downarrow & &\downarrow & &   \\
& & 0 &\rightarrow &S_3T^*\otimes T&\rightarrow &T^*\otimes F_0 &\rightarrow  & F_1& \rightarrow & 0 \\
& & \downarrow & & \downarrow  & & \downarrow & & \parallel  & & \\
0 & \rightarrow & R_3  & \rightarrow &J_3(T) &\rightarrow  & J_1(F_0)& \rightarrow &F_1 &\rightarrow & 0 \\
& & \downarrow& &\downarrow & &\downarrow & &\downarrow & &   \\
0 & \rightarrow & R_2 & \rightarrow & J_2(T)& \rightarrow & F_0 & \rightarrow & 0  & & \\
& & \downarrow & &\downarrow  & & \downarrow & & & &  \\
& & 0& & 0 & &  0  & & & &   
\end{array}     \]

\[  \begin{array}{rcccccccccl}
& & & & 0 & & 0 & & 0  & &  \\
& & & &\downarrow & &\downarrow & &\downarrow & &   \\
& & 0 & \rightarrow & 8  &\rightarrow &18 &\rightarrow  & 10 & \rightarrow & 0 \\
& & \downarrow & & \downarrow  & & \downarrow & & \parallel  & & \\
0 & \rightarrow & 3  & \rightarrow &  20 &\rightarrow  & 27 & \rightarrow & 10 &\rightarrow & 0 \\
& & \parallel  & &\downarrow & &\downarrow & &\downarrow & &   \\
0 & \rightarrow & 3 & \rightarrow &  12  & \rightarrow & 9 & \rightarrow & 0  & & \\
& & \downarrow & &\downarrow  & & \downarrow & & & &  \\
& & 0& & 0 & &  0  & & & &   
\end{array}     \]

\[  \hspace{15mm}SPENCER  \]
\[   \begin{array}{rcccccccccccl}
&&&&&&0& &  0 &       &  0 &   \\
&&&&&&\downarrow & & \downarrow  &  & \downarrow &  \\
&&0&\rightarrow&\Theta&\stackrel{j_q}{\longrightarrow}&C_0  &\stackrel{D_1}{\longrightarrow} & C_1 & \stackrel{D_2}{\longrightarrow} & C_2 & \rightarrow 0  \\
&&&&&&\downarrow & & \downarrow &   & \downarrow & \\
&&0 &\rightarrow  & T & \stackrel{j_q}{\longrightarrow} & C_0(T)&\stackrel{D_1}{\longrightarrow} & C_1(T) & \stackrel{D_2}{\longrightarrow} & C_2(T) & \rightarrow 0 \\
&&&&\parallel&& \hspace{5mm} \downarrow {\Phi}_0& & \hspace{5mm}\downarrow {\Phi}_1 &  & \hspace{5mm} \downarrow {\Phi}_2& \\
0 &\longrightarrow &\Theta &\rightarrow &T &\stackrel{{\cal{D}}}{\longrightarrow} &F_0  &\stackrel{{\cal{D}}_1}{\longrightarrow} & F_1 & \stackrel{{\cal{D}}_2}{\longrightarrow}  & F_2 & \rightarrow 0 \\
&&&&&& \downarrow  &  & \downarrow &  & \downarrow &    \\
& & & &&&0& & 0 &   & 0&
\end{array}  \]
\[  \hspace{17mm}JANET  \]

 \[  \begin{array}{rccccccccccccr}
  & 0& \longrightarrow& \tilde{\Theta} &\stackrel{j_2}{\longrightarrow}& 4 &\stackrel{D_1}{\longrightarrow}& 8 &\stackrel{D_2}{\longrightarrow} &  4 &\longrightarrow  0 &  \\

   & 0& \longrightarrow& \Theta &\stackrel{j_2}{\longrightarrow}& 3 &\stackrel{D_1}{\longrightarrow}& 6 &\stackrel{D_2}{\longrightarrow} &  3 &\longrightarrow  0 & \hspace{3mm}SPENCER  \\
  &&&&& \downarrow & & \downarrow & & \downarrow & &    & \\
   & 0 & \longrightarrow &   2   & \stackrel{j_2}{\longrightarrow} &  12  & \stackrel{D_1}{\longrightarrow} &   16  &\stackrel{D_2}{\longrightarrow} &  6   &   \longrightarrow 0 &   \hspace{5mm}  HYBRID\\
   & & & \parallel && \hspace{5mm}\downarrow {\Phi}_0 & &\hspace{5mm} \downarrow {\Phi}_1 & & \hspace{5mm}\downarrow {\Phi}_2 &  &\\
   0 \longrightarrow & \Theta &\longrightarrow &   2  & \stackrel{\cal{D}}{\longrightarrow} &  9  & \stackrel{{\cal{D}}_1}{\longrightarrow} & 10 & \stackrel{{\cal{D}}_2}{\longrightarrow} &   3  & \longrightarrow  0 & \hspace{7mm} JANET \\
   
    0 \longrightarrow & \tilde{\Theta} &\longrightarrow &   2  & \stackrel{{\cal{D}}}{\longrightarrow} &  8  & \stackrel{{\cal{D}}_1}{\longrightarrow} & 8  & \stackrel{{\cal{D}}_2}{\longrightarrow} &   2  & \longrightarrow  0 & 
   \end{array}     \]

\[   \begin{array}{rccccccccl}
& &  &    &  0  &  &  0 &  &    \\
 &   &  &   &   \downarrow  &  &  \downarrow  &  &    \\
  &  &  0 & \rightarrow  &   T^*\otimes ({\tilde{C}}_0/C_0) & = &  {\tilde{C}}_1/C_1 & \rightarrow & 0  \\
  &   &   \downarrow  &  &  \downarrow  &  &  \downarrow   &   & \\
0 & \rightarrow  & S_3T^*\otimes T & \rightarrow &T^*\otimes F_0 &  \rightarrow & F_1 & \rightarrow  & 0  \\
   &   &   \parallel &  &  \downarrow  &  &  \downarrow  &   & \\
0 & \rightarrow  & S_3T^*\otimes T& \rightarrow & T^*\otimes {\tilde{F}}_0 &\rightarrow & {\tilde{F}}_1 & \rightarrow  & 0 \\
  &   &   \downarrow  &  &  \downarrow  &  &  \downarrow   &   & \\
  &   &  0  &  &  0 &  & 0  &   & 
\end{array}   \]

\[   \begin{array}{rccccccccl}
&   &    &  &  0 &  & 0 &   &   \\
 &   &    &  &   \downarrow  &  &  \downarrow &  &   \\
&   & 0 & \rightarrow &  2 & \rightarrow & 2  & \rightarrow  & 0    \\
  &   &   \downarrow  &  &  \downarrow  &  &  \downarrow   &   & \\
0 & \rightarrow  & 8& \rightarrow &  18  &  \rightarrow & 10& \rightarrow  & 0  \\
   &   &   \parallel  &  &  \downarrow  &  &  \downarrow  &   & \\
0 & \rightarrow  & 8& \rightarrow & 16 &\rightarrow &8 & \rightarrow  & 0 \\
  &   &   \downarrow  &  &  \downarrow  &  &  \downarrow   &   & \\
  &   &  0  &  &  0 &  & 0  &   & 
\end{array}   \]

\[ \begin{array}{lcl}
{\tilde{C}}_r/C_r  \simeq  ker (F_r \rightarrow {\tilde{F}}_r)  \simeq {\wedge}^rT^*\otimes  ({\tilde{C}}_0/C_0) & \Rightarrow &
dim({\tilde{C}}_r) - dim(C_r)= dim (F_r) - dim ({\tilde{F}}_r)   \\
     & \Rightarrow &  dim (C_r) + dim (F_r) = dim ({\tilde{C}}_r) + dim ({\tilde{F}}_r)
     \end{array}    \]

In this new situation, we now notice that ${\tilde{R}}_2\subsetneq {\rho}_1({\tilde{R}}_1)=J_1({\tilde{R}}_1)\cap J_2(T)\subset J_1(J_1(T))$ and the induced morphism ${\tilde{F}}_0\rightarrow J_1({{\tilde{F'}}}_0)$ is thus no longer a monomorphism though we still have an isomorphism ${\tilde{F}}"_0 \simeq S_2T^*\otimes T$ because ${\tilde{g}}_2=0$ again. Finally, we may extend such a procedure to the conformal group of space-time
by considering the system of infinitesimal conformal transformations of the Minkowski metric defined by the first order system ${\hat{R}}_1 \subset J_1(T)$ in such a way that we have the strict inclusions $R_2 \subset {\tilde{R}}_2 \subset {\hat{R}}_2 \subset J_2(T)$ with $dim({\hat{g}}_2)=n=4$. For this, we just need to introduce the metric density $ \hat{\omega}= \omega (\mid det(\omega)\mid)^{ - \frac{1}{n}}$ and consider the system 
${\cal{L}}(\xi)\hat{\omega}=0$ ([38]): \\
\[   {\hat{\omega}}_{rj}{\partial}_i{\xi}^r + {\hat{\omega}}_{ir}{\partial}_j{\xi}^r  - \frac{2}{n}{\hat{\omega}}_{ij}{\partial}_r{\xi}^r + {\xi}^r {\partial}_r  {\hat{\omega}}_{ij} = 0        \]

\noindent
{\bf A) CONTROL THEORY}:  \\
\noindent
{\bf EXAMPLE 4.A.1}: ({\it OD control theory}) In classical control theory we have $n=1$ and the only independent variable is the time, simply denoted by $x$ but we may choose any ground differential field like $K=\mathbb{Q}(x)$. In that case, we shall refer to ([22] or [25]) for the proof of the following technical results that will be used in this case (Compare to [49]). Instead of the standard "upper dot" notation for derivative we shall identify the formal and the jet notations, setting thus $d_xy=dy=y_x$. With $m=2$, let us consider the elementary Single Input/Single Output (SISO) second order system $y^1_{xx} - y^2_x+a(x)y^2=0$ with a variable coefficient $a\in K$. The corresponding formally surjective operator is ${\partial}_{xx}{\eta}^1- {\partial}_x{\eta}^2 + a(x){\eta}^2=\zeta$. Treating such a system by using classical methods is not so easy when $a$ is not constant as it cannot be possible to transform it to the standard Kalman form. On the contrary, multiplying by a test function (or Lagrange multiplier) $\lambda$ and integrating by parts, we obtain the adjoint system/operator:  \\
\[ \left  \{  \begin{array}{lclcl}
y^1 & \longrightarrow & {\lambda}_{xx}& = & {\mu}^1  \\
y^2 & {\longrightarrow}&  {\lambda}_x +a \lambda& = &{\mu}^2
\end{array}  \right.  \fbox{ $ \begin{array}{ll}
1   \\
 \bullet 
\end{array} $ }  \]
This system has a trivially involutive zero symbol but is not even formally integrable and we have to consider :  \\
\[ \left  \{  \begin{array}{lcl}
 {\lambda}_{xx}  & = & {\mu}^1  \\
  {\lambda}_x +a \lambda  & = &{\mu}^2  \\
  ({\partial}_xa - a^2)\lambda & = & {\mu}^2_x - {\mu}^1- a {\mu}^2
\end{array}  \right.  \fbox{ $ \begin{array}{ll}
1   \\
 \bullet \\
 \bullet
\end{array} $ }  \]
We have thus two possibilities:  \\
\noindent
$\bullet$ We have $a_x-a^2\neq 0$ and the adjoint system has the only zero solution, that is the adjoint operator is injective. In this case $N=0$ and thus $t(M)=ext^1(N)=0$ that is $M$ is torsion-free. However, as $n=1$ it follows that $D=K[d]$ is a principal ideal ring which is therefore free and thus projective ([Kunz,Rot]), that is $M$ {\it is torsion-free if and only if} $N=0$ and the system is controllable.  \\

\noindent
$\bullet$  The Riccati equation $a_x- a^2=0$ is satisfied, for example if $a= -1/x$ and we get the CC ${\mu}^2_x - {\mu}^1 - a {\mu}^2=0$. Multiplying by a test function $\xi$ and integrating by parts, we get the adjoint operator:    \\
\[ \left  \{  \begin{array}{lclcl}
{\mu}^1 & \longrightarrow & - {\xi} & = & {\eta}^1  \\
{\mu}^2 & {\longrightarrow}& - {\xi}_x -  a \xi  & = &{\eta}^2
\end{array}  \right. \]
with only one first order generating CC, namely ${\partial}_x{\eta}^1 - {\eta}^2 +a {\eta}^1=0$. It follows that $N\neq 0 \Rightarrow ext^1(N)\neq 0$ is a torsion module generated by the residue of $z=y^1_x - y^2 + a y^1$. We obtain indeed a torsion element as we can check at once that $z_x - az=0$ and wish good luck for control people to recover this result even on such an elementary example because the Kalman criterion is only working for systems with constant coeficients (Compare [25] and [49]).  \\

\noindent
{\bf EXAMPLE 4.A.2}: ({\it PD control theory}) With $n=2$, let us consider the (trivially involutive) inhomogeneous single first order PD equations with two independent variables $(x^1,x^2)$, two unknown functions $({\eta}^1, {\eta}^2)$ and a second member $\zeta$: \\
\[ {\partial}_2{\eta}^1 - {\partial}_1{\eta}^2 + x^2 {\eta}^2=\zeta   \hspace{1cm} \Leftrightarrow \hspace{1cm}  {\cal{D}}_1 \eta = \zeta \]
The ring of differential operators is $D=K[d_1,d_2]$ with $K=\mathbb{Q}(x^1,x^2)$. Multiplying on the left by a test function $\lambda$ and integrating by parts, the corresponding adjoint operator is described by:  \\
\[  \left\{ \begin{array}{rcll}
{\eta}^1 & \rightarrow  & - {\partial}_2 \lambda  & ={\mu}^1  \\
{\eta}^2 & \rightarrow & \hspace{3mm}{\partial}_1 \lambda + x^2 \lambda  & ={\mu}^2
\end{array} \right.  \hspace{1cm} \Leftrightarrow \hspace{1cm}  ad({\cal{D}}_1)\lambda=\mu\]
Using crossed derivatives, this operator is injective because $\lambda={\partial}_2{\mu}^2+{\partial}_1{\mu}^1+x^2{\mu}^1$ and we even obtain a lift  $\lambda \longrightarrow \mu \longrightarrow \lambda$. Substituting, we get the two CC:  \\
\[  \left  \{  \begin{array}{lcl}
  {\partial}_{22}{\mu}^2 + {\partial}_{12}{\mu}^1+x^2{\partial}_2{\mu}^1 + 2{\mu}^1  & = & {\nu}^1\\
{\partial}_{12}{\mu}^2 + {\partial}_{11}{\mu}^1+2x^2{\partial}_1{\mu}^1+x^2{\partial}_2{\mu}^2+(x^2)^2{\mu}^1-{\mu}^2 & = & {\nu}^2  
\end{array}  \right.  \fbox{ $ \begin{array}{ll}
1 & 2   \\
1 & \bullet 
\end{array} $ } \]
This system is involutive and the corresponding generating CC for the second member $({\nu}^1,{\nu}^2)$ is:  \\
\[  {\partial}_2{\nu}^2 - {\partial}_1{\nu}^1 - x^2 {\nu}^1=0  \]
Therefore ${\nu}^2$ is differentially dependent on ${\nu}^1$ but ${\nu}^1$ is also differentially dependent on ${\nu}^2$. Multiplying on the left by a test function $\theta$ and integrating by parts, the corresponding adjoint system of PD equations is:  \\
\[  \left\{ \begin{array}{rcll}
{\nu}^1 & \rightarrow  & \hspace{3mm}{\partial}_1 \theta - x^2 \theta & ={\xi}^1  \\
{\nu}^2 & \rightarrow & -{\partial}_2 \theta   & ={\xi}^2
\end{array} \right.  \hspace{1cm} \Leftrightarrow \hspace{1cm}  ad({\cal{D}}_{-1})\theta=\xi  \]
Multiplying now the first equation by the test function ${\xi}^1$, the second equation by the test function ${\xi}^2$, adding and integrating by parts, we get the {\it canonical parametrization} ${\cal{D}}\xi = \eta $:  \\
\[   \left \{\begin{array}{rcll}
{\mu}^2 & \rightarrow & {\partial}_{22}{\xi}^1+{\partial}_{12}{\xi}^2-x^2{\partial}_2{\xi}^2-2{\xi}^2  & = {\eta}^2  \\
{\mu}^1&  \rightarrow & {\partial}_{12}{\xi}^1-x^2{\partial}_2{\xi}^1+{\xi}^1+{\partial}_{11}{\xi}^2  -2x^2{\partial}_1{\xi}^2+(x^2)^2{\xi}^2 & = {\eta}^1        
\end{array}  \right. \fbox{ $ \begin{array}{ll}
1 & 2   \\
1 & \bullet 
\end{array} $ } \]
of the initial system with zero second member. This system is involutive and the kernel of this parametrization has differential rank equal to $1$ because ${\xi}^1$ or ${\xi}^2$ can be given arbitrarily.  \\
Keeping now ${\xi}^1=\xi$ while setting ${\xi}^2=0$, we get the {\it first second order minimal parametrization} $\xi \rightarrow ({\eta}^1,{\eta}^2)$:\\
\[  \left  \{  \begin{array}{ll}
{\partial}_{22}\xi & = {\eta}^2  \\
{\partial}_{12}\xi-x^2{\partial}_2\xi+\xi & ={\eta}^1
\end{array}  \right.  \fbox{ $ \begin{array}{ll}
1 & 2   \\
1 & \bullet 
\end{array} $ }  \]
This system is again involutive and the parametrization is minimal because the kernel of this parametrization has differential rank equal to $0$. With a similar comment, setting now ${\xi}^1=0$ while keeping ${\xi}^2={\xi}'$, we get the {\it second second order minimal parametrization} ${\xi}' \rightarrow ({\eta}^1,{\eta}^2)$:  \\
\[  \left  \{\begin{array}{ll}
{\partial}_{11}{\xi}'-2x^2{\partial}_1{\xi}'+(x^2)^2{\xi}' & = {\eta}^1  \\
{\partial}_{12}{\xi}'  - x^2{\partial}_2{\xi}' - 2{\xi}' & = {\eta}^2 
\end{array}   \right.  \]
which is again easily seen to be involutive by exchanging $x^1$ with $x^2$.  \\
With again a similar comment, setting now ${\xi}^1={\partial}_1\phi, {\xi}_2= - {\partial}_2\phi$ in the canonical parametrization, we obtain the {\it third different second order minimal parametrization}:  \\
\[  \left  \{  \begin{array}{ll}
x^2{\partial}_{22}\phi + 2 {\partial}_2\phi& = {\eta}^2  \\
x^2{\partial}_{12}\phi - (x^2)^2{\partial}_2\phi + {\partial}_1\phi& ={\eta}^1
\end{array}  \right.  \fbox{ $ \begin{array}{ll}
1 & 2   \\
1 & \bullet 
\end{array} $ }  \]
We are now ready for understanding the meaning and usefulness of what we have called " {\it relative parametrization} " in ([29]) by imposing the {\it differential constraint} ${\partial}_2{\xi}^1 + {\partial}_1{\xi}^2=0$ which is compatible as we obtain indeed the new {\it first order relative parametrization}:  \\
\[  \left  \{  \begin{array}{ll}
{\partial}_2{\xi}^1 +{\partial}_1{\xi}^2 &  =  0  \\
-x^2{\partial}_2{\xi}^2 - 2 {\xi}^2  & = {\eta}^2  \\
-x^2{\partial}_1{\xi}^2 + (x^2)^2{\xi}^2 + {\xi}^1  & ={\eta}^1
\end{array}  \right.  \fbox{ $ \begin{array}{ll}
1 & 2   \\
1  &  2  \\
1 & \bullet 
\end{array} $ }  \]
with  $2$ equations of class $2$ (thus with class $2$ full) and only $1$ equtaion of class $1$.  \\
In a different way, we may add the differential constraint ${\partial}_1{\xi}^1+{\partial}_2{\xi}^2=0$ but we have to check that it is compatible with the previous parametrization. For this, we have to consider the following second order system which is easily see to be involutive with $2$ second order equations of (full) class $2$, (only) $2$ second order equations of class $1$ and $1$ equation of order $1$:  \\
\[   \left \{\begin{array}{lcl}
{\partial}_{22}{\xi}^2 + {\partial}_{12}{\xi}^1 & = & 0  \\
{\partial}_{22}{\xi}^1+{\partial}_{12}{\xi}^2-x^2{\partial}_2{\xi}^2-2{\xi}^2  & = &{\eta}^2  \\
{\partial}_{12}{\xi}^2 + {\partial}_{11}{\xi}^1 & = & 0  \\
 {\partial}_{12}{\xi}^1-x^2{\partial}_2{\xi}^1+{\xi}^1+{\partial}_{11}{\xi}^2  -2x^2{\partial}_1{\xi}^2+(x^2)^2{\xi}^2 & = &{\eta}^1   \\
{\partial}_2{\xi}^2 +{\partial}_1 {\xi}^1   &  = & 0 
\end{array}  \right. \fbox{ $ \begin{array}{ll}
1 & 2   \\
1 & 2 \\
1 & \bullet \\
1 & \bullet\\
\bullet & \bullet  
\end{array} $ } \]
The $4$ generating CC only produce the desired system for $({\eta}^1,{\eta}^2)$ as we wished.  \\
We cannot impose the condition ${\cal{D}}_{-1}\theta=\xi $ already found as it should give the identity $0 = \eta$.  \\
It is however also important to notice that the strictly exact long exact sequence:
\[   0  \longrightarrow D \stackrel{{\cal{D}}_1}{\longrightarrow} D^2 \stackrel{{\cal{D}}}{\longrightarrow} D^2 \stackrel{{\cal{D}}_{-1}}{\longrightarrow} D \longrightarrow 0  \]
splits because we have a lift $\zeta \longrightarrow \eta \longrightarrow \zeta$, namely:   \\
\[ \zeta \longrightarrow ( -{\partial}_1 \zeta  + x^2\zeta = {\eta}^1, - {\partial}_2 {\eta}^2 ={\eta}^2) \longrightarrow {\partial}_2{\eta}^1 - {\partial}_1{\eta}^2 + x^2 {\eta}^2=\zeta   \] 
We have thus an isomorphism $D^2 \simeq D \oplus M$ in the resolution $0 \longrightarrow D  \stackrel{{\cal{D}}_1}{\longrightarrow} D^2 \stackrel{p}{\longrightarrow} M  \longrightarrow 0 $ and all the differential modules defined from the operators involved are projective, thus torsion-free or $0$-pure with vanishing $r$-extension modules $ext^r( \,\,\, )=0,\forall r\geq 1$.   \\
As an exercise, we finally invite the reader to study the situation met with the system ${\partial}_2{\eta}^1 - {\partial}_1{\eta}^2 + a(x){\eta}^2$ whenever $a\in K$ ({\it Hint}: The controllability condition is now ${\partial}_1a\neq 0$). The comparison with the previous OD case needs 
no comment. \\ 

\noindent
{\bf B) ELECTROMAGNETISM}: \\
Most physicists know the Maxwell equations in vacuum, eventually in dielectrics and magnets, but are largely unaware of the more delicate constitutive laws involved in field-matter couplings like piezzoelectricity, photoelasticity or streaming birefringence. In particular they do not know that the phenomenological laws of these phenomena have been given ... by Maxwell ([37]). The situation is even more critical when they deal with invariance properties of Maxwell equations because of the previous comments ([5]). Therefore, we shall first quickly recall what the use of adjoint operators and differential duality can bring when studying Maxwell equations as a first step before providing comments on the so-called gauge condition brought by the danish physicist Ludwig Lorenz in 1867 and not by Hendrik Lorentz with name associated with the Lorentz transformations.  \\
Theough it is quite useful in actual practice, the following approach to Maxwell equations cannot be found in any textbook. Namely, avoiding any variational calculus based on given Minkowski constitutive laws ${\cal{F}}\sim F$ between {\it field} $F$ and {\it induction} ${\cal{F}}$ for dielectric or magnets, let us  use differential duality and define the first set $M_1$ of Maxwell equations by $d:{\wedge}^2T^* \rightarrow {\wedge}^3T^*$ while the second set $M_2$ will be defined by $ad(d):{\wedge}^4T^*\otimes {\wedge}^2T \rightarrow {\wedge}^4T^*\otimes T$ with $d:T^*\rightarrow {\wedge}^2T^*$, in a totally independent and intrinsic manner, {\it using now contravariant tensor densities in place of covariant tensors}. As we have already proved since a long time in ([20-22],[24],[37-38]), the key result is that {\it these two sets of Maxwell equations are invariant by any diffeomorphism}, contrary to what is generally believed ([5]). We recapitulate below this procedure in the form of a (locally exact) differential sequence and its (locally exact) formal adjoint sequences where the left dotted arrow is the standard composition of operators:  \\
\[  \begin{array}{ccccccc}
potential &\stackrel{d}{\longrightarrow}& field& \rightarrow & induction &\stackrel{ad(d)}{\longrightarrow}& current  \\
A &\longrightarrow &F& \rightarrow &{\cal{F}}& \longrightarrow &{\cal{J}}
\end{array}   \]
which is responsible for EM waves, though it is equivalent  to the composition:   \\
\[   pseudopotential \stackrel{ad(d)}{\longrightarrow} induction \rightarrow field \stackrel{d}{\longrightarrow} {\wedge}^3T^*  \]
The main difference is that we need to set ${\cal{J}}=0$ in the first approach because of $M_2$ while we get {\it automatically} such a vanishing  assumption in the second approach because of $M_1$, avoiding therefore the Lorenz condition as in ([35], Remark 5.5). \\
\[   \begin{array}{ccccc}
  potential= (A_i) &\stackrel{d}{\longrightarrow}& ({\partial}_iA_j - {\partial}_jA_i=F_{ij})=field &\stackrel{M_1}{\longrightarrow}& ({\partial}_iF_{jk} + {\partial}_jF_{ki} + {\partial}_kF_{ij}=0)   \hspace{1cm}     \\
   &  &  &  &  \\
   current=({\partial}_i{\cal{F}}^{ij}={\cal{J}}^j )& \stackrel{M_2}{\longleftarrow}& ({\cal{F}}^{ij}) =induction & \stackrel{ad(d)}{\longleftarrow } &  pseudopotential  
\end{array}   \]

\[   \begin{array}{ccccc}
A  &  &  F  &  &  \\
T^* & \stackrel{d}{\longrightarrow} & {\wedge}^2T^*  &  \stackrel{d=M_1}{\longrightarrow} & {\wedge}^3T^*  \\
  \vdots &   &\downarrow   &   & \vdots  \\
{\wedge}^4T^*\otimes T & \stackrel{ad(d)=M_2}{\longleftarrow} & {\wedge}^4T^*\otimes {\wedge}^2T&  \stackrel{ad(d)}{\longleftarrow}& {\wedge}^4T^*\otimes {\wedge}^3T   \\
 \updownarrow &  & \updownarrow &  & \updownarrow \\
  {\wedge}^3T^* &  \stackrel{d}{\longleftarrow} & {\wedge}^2T^*&\stackrel{d}{\longleftarrow}& T^*  \\
  {\cal{J}} &&  {\cal{F}}  &  &  
\end{array}    \]

Using symbolic notations with an euclidian metric instead of the Minkowski one because they are both locally constant while using the constitutive law ${\cal{F}}=F$ for simplicity in vacuum while raising or lowering the indices by means of the metric, we have the parametrization $d_iA_j - d_jA_i=F_{ij}$ and obtain by composition in the left upper square:   \\
\[ d_i(d^iA^j - d^jA^i)=d_id^iA^j - d^j(d_iA ^i)= {\cal{J}}^j \Rightarrow d_j(d_id^iA^j - d^jd_iA^i)= d_j{\cal{J}}^j=0   \]
with implicit summations on $i$ and $j$. We may consider the composit homogeneous second order system $d_{ii}A_j - d_{ij}A_i=0$ which is automatically formally integrable and is easily seen (exercise) to be involutive. The character ${\alpha}^n_2$ is obtained by considering 
$d_{nn}A_j - d_{jn}A_n$ for the equation giving ${\cal{J}}_j$. For ${\cal{J}}_n$ we get $d_{nn}A_n - d_{nn}A_n=0$ and thus ${\alpha}^n_2=n-(n-1)=1$ a result showing that the corresponding differential module has rank $1$ and there is thus only one CC, namely $d_j{\cal{J}}^j=0$ with implicit summation on $j$. We prove that we may add the Lorenz condition $d_iA^i=0$ to bring the rank to zero. Indeed, we have now the inhomogeneous system $d_id^iA^j={\cal{J}}^j$ and the differential constraint thus brought is compatible with the conservation of current. The corresponding homogeneous system obtained by adding the Lorenz constraint has second order symbol obtained by considering {\it both} $d_id^iA^j=0$ {\it and} $ d_{ij}A ^i=0$ 
{\it or} $d_{ij}A^j=0$. We obtain therefore $d_{nn}A_j=0, d_{nn}A_n=0$ showing that we have now ${\alpha}^n_2=0$ and a torsion differential module. As a more important and effective result that does not seem to be known, we have:  \\
\noindent
{\bf PROPOSITION 4.B.1}: When $n=4$, the system:  \\
\[   \left \{\begin{array}{lcccl}
{ \Psi }^j & \equiv &  d_{44}A^j + ... +d_{11}A^j & = & {\cal{J}}^j  \\
  d_4 \Phi - {\Psi}^4 &  \equiv & d_{34}A^3 +d_{24}A^2 + d_{14} A^1 - d_{33}A^4 - d_{22}A^2 - d_{11}A^4 & = & - {\cal{J}}^4 \\
  d_3 \Phi & \equiv &  d_{34}A^4 + ... + d_{13}A^1& = & 0  \\
 d_2 \Phi & \equiv &   d_{24}A^4 + ... + d_{12}A^1& = & 0 \\
 d_1\Phi & \equiv &  d_{14}A^4 + ... + d_{11}A^1 & = &  0  \\
 \Phi & \equiv &   d_4A^4 + ...  + d_1A^1&  = & 0
\end{array}  \right. \fbox{ $ \begin{array}{cccc}
1 & 2 & 3 & 4  \\
1 & 2 & 3 & \bullet \\
1 & 2 & 3 & \bullet  \\
1 & 2 & \bullet & \bullet \\
1  &\bullet & \bullet & \bullet\\
\bullet & \bullet &\bullet & \bullet 
\end{array} $}  \]
is involutive with four equations of class $4$, two equations of class 3, one equation of class $2$ and one equation of class $1$. The $11$ resulting CC only provide the conservation of current.  \\

\noindent
{\it Proof}: Using the corresponding {\it Janet tabular} on the let, one can check at once that the $4$ CC brought by the only first order equation $\Phi=0$ do not bring anything new, as they amount to crossed derivatives, and that {\it we are only left with the $4$ upper dots on the right side}. However, for $i=1,2,3$, we have $d_4(d_i\Phi)=d_i(d_4\Phi - {\Psi}^4) + d_i {\Psi}^4$ and we are thus only left with a single CC, getting successively:  \\
\[  \begin{array}{lcl}
d_4 (d_4\Phi - {\Psi}^4) & \equiv &  d_{344} A^3 + d_{244}A^2 + d_{144}A^1 - d_{334}A^4 - d_{224}A^2 - d_{114}A^1   \\
- d_3 {\Psi}^3 & \equiv &  -d_{344}A^3 - d_{333}A^3 - d_{223}A^3 - d_{113} A^3   \\
-d_2{\Psi}^2 & \equiv & -d_{244}A^2 - d_{233} A^2 - d_{222}A^2 - d_{112}A^2    \\
- d_1{\Psi}^1 & \equiv & -d_{144} A^1 - d_{133} A^1 - d_{122}A^1 - d_{111}A^1   \\
 d_3(d_3\Phi) & \equiv &  d_{334}A^4 + d_{333}A^3 + d_{233} A^2 + d_{133}A^1  \\
 d_2(d_2 \Phi) & \equiv &  d_{224}A^4 + d_{223}A^3 + d_{222}A^2 + d_{122} A^1 \\
 d_1(d_1\Phi) & \equiv &   d_{114}A^4 + d_{113}A^3 +d_{112}A^2 +d_{111}A^1
\end{array}  \]
Summing these $7$ equations, we are left with the identity $ -(d_{44} \Phi + ...  +d_{11}\Phi ) + d_j{\Psi}^j=d_j{\cal{J}}^j=0$. It is important to notice that {\it no other procedure} can prove that we have an involutive symbol in $\delta$-regular coordinates and this is the only way to compute {\it effectively} all the four characters $(0< 6 < 11 < 15)$ with $6+11+15=32= (4\times 10) -(4+4)$ for the dimension of the symbol of order $2$, a result not evident at first sight. Accordingly, the so-called Lorenz gauge condition is only a pure " {\it artifact} " amounting to a relative minimum parametrization with no important physical meaning as it can be avoided by using only the EM field $F$ ([31],[34],[35]).    \\
\hspace*{12cm} Q.E.D.   \\
Such a new approach to a classical result is nevertheless bringing a totally unsatisfactory consequence. Using the well known correspondence between electromagnetism(EM) and elasticity (EL) used for all engineering computations with finite elements:  \\
\[EM \,\,potential \leftrightarrow EL\,\, displacement, \hspace{5mm}  EM \,\,field \leftrightarrow EL \,\,strain, \hspace{5mm} EM \,\,induction \leftrightarrow EL \,\,stress \]
where EL means elasticity, and instead of the left upper square in the diagram, ... {\it we have to consider the right upper square}.  \\
We finally prove that the use of the linear {\it and} nonlinear Spencer operators drastically changes the previous standard procedure in a way that could not even be imagined with classical methods. For such a purpose, we make a few comments on the implicit summation appearing in differential duality. For example, we have, up to a divergence: \\
\[  {\cal{X}}^{,r}_kX^k_{,r}= {\cal{X}}^{,r}_k( {\partial}_r{\xi}^k- {\xi}^k_r)= - {\partial}_r({\cal{X}}^{,r}_k){\xi}^k - {\cal{X}}^{,r}_k{\xi}^k_r + ... \] 
In the conformal situation, we have ${\xi}^1_1={\xi}^2_2= ... ={\xi}^n_n= \frac{1}{n}{\xi}^r_r$ and obtain therefore, as factor of the firs jets:  \\
\[  {\cal{X}}^{,1}_1 {\xi}^1_1+ {\cal{X}}^{,2}_2 {\xi}^2_2+ ... + {\cal{X}}^{,n}_n{\xi}^n_n = ({\cal{X}}^{,r}_r)\frac{1}{n}{\xi}^r_r= ({\cal{X}}^{,r}_r){\xi}^1_1  \]
Going to the next order, we get as in ([26]), up to a divergence:  \\
\[  {\cal{X}}^{1,r}_1{\partial}_r{\xi}^1_1= - ({\partial}_r{\cal{X}}^{1,r}_1){\xi}^1_1 +  ...  \]
Collecting the results and changing the sign, we obtain for the first time the {\it Cosserat equation for the dilatation}, namely the so-called {\it virial equation} that we provided in $2016$ ([34], p 35) :  \\
\[   {\partial}_r {\cal{X}}^{1,r}_1 + {\cal{X}}^{,r}_r= 0       \]
generalizing the well known {\it Cosserat equations for the rotations} provided in $1909$ ([6], p 137):  \\
\[             {\partial}_r{\cal{X}}^{ij,r} + {\cal{X}}^{i,j} - {\cal{X}}^{j,i}=0   \]
As for EM, substituting ${\partial}_i{\xi}^r_{rj} - {\partial}_j{\xi}^r_{ri}$ in the dual sum ${\cal{F}}^{i<j}F_{i<j}=
 \frac{1}{2} {\cal{F}}^{ij}F_{ij}$ and integrating by parts, we get a part of the {\it Cosserat equations for the elations}, namely:ÊÊ\\
 \[      {\partial }_r{\cal{F}}^{ir} - {\cal{J}}^i=0  \Rightarrow {\partial}_i{\cal{J}}^i=0 \Rightarrow {\cal{X}}^{,r}_r=0     \]
saying that {\it the trace of the EM impulsion-energy tensor must vanish} ([34], p 37).  \\
We sum up all these results in the following tabular only depending on the Spencer operator:   \\

\[  \renewcommand \arraystretch{1.5}
\begin{tabular}{|c|c|c|}
\hline  
\multicolumn{2}{|c|}{FIELD}  & INDUCTION  \\  [3pt]
\hline
 NONLINEAR & LINEAR & DUAL \\  [3pt]
\hline  
$\bar{D}f_{q+1}={\chi}_q \in T^*\otimes R_q$ &  $ D{\xi}_{q+1}=X_q\in T^*\otimes R_q $ & $ {\cal{X}}_q \in {\wedge}^nT^*\otimes T \otimes R_q^*$\\ [1mm]
\hline
${\chi}^k_{,r} $   &$  {\partial}_r{\xi}^k -{\xi}^k_r=X^k_{,r} $                  & $ {\cal{X}}^{,r}_k  $  \\ [1mm]
${\chi}^k_{i,r} $  &$  {\partial}_r{\xi}^k_i - {\xi}^k_{i,r}=X^k_{i,r} $        & $ {\cal{X}}^{i,r}_k  $ \\ [1mm]
$ {\chi}^r_{r,i}$ & $ {\partial}_i{\xi}^r_r - {\xi}^r_{ri}= X^r_{r,i}=X_i $ & ${\cal{J}}^i $ \\  [1mm]
${\chi}^k_{ij,r} $ &$  {\partial}_r{\xi}^k_{ij} - {\xi}^k_{ijr}=  {\partial}_r{\xi}^k_{ij}=X^k_{ij,r} $ & $ {\cal{X}}^{ij,r}_k $ \\ [1mm]
${\partial}_i{\chi}^r_{r,j} - {\partial}_j{\chi}^r_{r,i}={\varphi}_{ij}$ & $ {\partial}_i{\xi}^r_{r,j} - {\partial}_j{\xi}^r_{r,i}=F_{ij}$ & 
${\cal{F}}^{ij} $  \\[1mm]
            & $ \frac{1}{2}({\partial}_i{\xi}^r_{r,j} + {\partial}_j{\xi}^r_{r,i})=R_{ij}$ &        \\[1mm]
\hline
\end{tabular}     \]    
We notice that $F=(F_{ij}) \in {\wedge}^2T^*, R=(R_{ij})\in S_2T ^*$ and the part of the first Spencer bundle made by the $1$-forms with value in the $n=4$ elations provides the splitting:   \\
 \[         (F,R) \in {\wedge}^2T^* \oplus S_2T^*\simeq T^*\otimes T^*  \]
because of the well known exactness of the Spencer $\delta$-sequence:  \\
\[  0 \longrightarrow S_2T^* \stackrel{\delta}{\longrightarrow} T^*\otimes T^* \stackrel{\delta}{\longrightarrow} {\wedge}^2T^* \longrightarrow 0 \]

\noindent
{\bf C) GENERAL RELATIVITY}:   \\
Roughly speaking, we shall say that a parametrization of an operator is {\it minimal} if its corresponding operator defines a torsion module or, equivalently, if the kernel of the parametrizing operator has differential rank equal to $0$. It is not so well known even today that, up to an isomorphism, the {\it Cauchy} stress operator essentially admits only one parametrization in dimension $n=2$ which is minimum but the situation is quite different in dimension $n=3$. Indeed, the parametrization found by E. Beltrami in $1892$ with $6$ potentials ([33]) is not minimal as the kernel of the {\it Beltrami} operator has differential rank $3$ while the {\it two other} parametrizations respectively found by J.C. Maxwell in $1870$ and by G. Morera in $1892$  are both minimal with only $3$ potentials even though they are quite different because the first is cancelling $3$ among the $6$ potentials while the other is cancelling the $3$ others. In particular, we point out the technical fact that it is quite difficult to prove that the Morera parametrization is providing an involutive system. These three tricky examples are proving that the possibility to exhibit different parametrizations of the stress equations that we have presented has surely nothing to do with the proper mathematical background of elasticity theory as it provides an explicit application of double differential duality in differential homological algebra. Also, the example presented in Section 3.A is proving that the existence of many different minimal parametrizations has surely nothing to do with the mathematical foundations of control theory. Similarly, we have just seen in the previous section that the so-called Lorenz condition has surely nothing to do with the mathematical foundations of EM. Such a comment will be now extended in a natural manner to GR. \\
With tandard notations, denoting by $\Omega \in S_2T^*$ a perturbation of the non-degenerate metric $\omega$, it is well known (See [30] and [35] for more details) that the linearization of the {\it Ricci} tensor $R=(R_{ij})\in S_2T^*$ over the Minkowski metric, considered as a second order operator $\Omega \rightarrow R$, may be written {\it with four terms} as:  \\
\[  2 R_{ij}= {\omega}^{rs}(d_{ij}{\Omega}_{rs}+d_{rs}{\Omega}_{ij}-d_{ri}{\Omega}_{sj} - d_{sj}{\Omega}_{ri})= 2R_{ji}  \]
Multiplying by test functions $({\lambda}^{ij}) \in {\wedge}^4T^*\otimes S_2T$ and integrating by parts on space-time, we obtain the following {\it four terms} describing the so-called {\it gravitational waves equations}:  \\
\[   (\Box {\lambda}^{rs}+ {\omega}^{rs}d_{ij}{\lambda}^{ij}-{\omega}^{sj}d_{ij}{\lambda}^{ri}- {\omega}^{ri}d_{ij}{\lambda}^{sj}){\Omega}_{rs}  =  {\sigma}^{rs}{\Omega}_{rs}\]
where $\Box$ is the standard Dalembertian. Accordingly, we have:   \\
\[   d_r{\sigma}^{rs}={\omega}^{ij}d_{rij} {\lambda}^{rs}+{\omega}^{rs}d_{rij}{\lambda}^{ij}-
{\omega}^{sj}d_{rij}{\lambda}^{ri} - {\omega}^{ri} d_{rij}{\lambda}^{sj}=0  \]
The basic idea used in GR has been to simplify these equations by adding the {\it differential constraints} $d_r{\lambda}^{rs}=0$ in order to find only 
$\Box {\lambda}^{rs}={\sigma}^{rs}$, exactly like in the Lorenz condition for EM. Before going ahead, it is important to notice that when $n=2$, the Lagrange multiplier $\lambda$ is just the Airy function $\phi$ and, using an integration by parts,  we have the identity:\\
\[   \phi (d_{11}{\Omega}_{22} +d_{22}{\Omega}_{11} - 2 d_{12}{\Omega}_{12})= d_{22}\phi {\Omega}_{11} - 2 d_{12}\phi {\Omega}_{12} + d_{11}{\Omega}_{22} + div( \,\,\,)  \]
providing the {\it Airy parametrization} of the Cauchy stress equations: \\
\[   {\sigma}^{11}=d_{22}\phi, \,\, \,\,{\sigma}^{12}={\sigma}^{21}= - d_{12}\phi, \,\,\,\, {\sigma}^{22}=d_{11}\phi  \]
where {\it the Airy function has nothing to do with the perturbation of the metric}.  \\
However, even if its clear that the constraints are compatible with the Cauchy equations, we do believe that the following result is not known as it does not contain any reference to the usual {\it Einstein} tensor $E_{ij}=R_{ij}- \frac{1}{2}{\omega}_{ij}tr(R)$ where $tr(R)={\omega}^{rs}R_{rs}$, which is therefore {\it useless} because it contains $6$ terms instead of $4$ terms only.  \\

\noindent
{\bf PROPOSITION 4.C.1}:  The system made by $\Box {\lambda}^{rs}={\sigma}^{rs}$ and $d_r{\lambda}^{rs}=0$ is a relative minimum involutive parametrization of the Cauchy equations describing the formal adjoint of the Killing operator, that is $Cauchy=ad(Killing)$ as operators.  \\

\noindent
{\it Proof}: For each {\it given} $s=1,2,3,4$ the system under study is {\it exactly} the system used for studying the Lorenz condition in Proposition 4.B.1. Accordingly, nothing has to be changed in the proof of this proposition and we get an involutive second order sysem with $d_r{\sigma}^{rs}=0$ as only CC in place of the conservation of current. Needless to say that this result has nothing to do with any concept of gauge theory as it is sometimes claimed ([8],[30]).  \\
\hspace*{12cm}   Q.E.D.   \\

\noindent
{\bf 5) CONCLUSION}:  \\
In 1916, F.S. Macaulay used a new localization technique for studying unmixed polynomial ideals. In 2013, we have generalized this procedure in order to study pure differential modules, obtaining therefore a relative parametrization in place of the absolute parametrization already known for torsion-free modules and equivalent to controllability in the study of OD or PD control systems, a result showing that {\it controllability does not depend on the choice of the control variables}, despite what engineers still believe. Meanwhile, we have pointed out the existence of minimum parametrizations obtained by adding in a convenient but generally not intrinsic way certain compatible differential constraints on the potentials. We have proved that this is {\it exactly} the kind of situation met in control theory, in EM with the Lorenz condition and in GR with gravitational waves. However, the systematic use of {\it adjoint operators} and {\it differential duality} is proving that the physical meaning of the potentials involved has {\it absolutely nothing to do} with 
the one usually adopted in these domains. Therefore, these results bring the need to revisit the mathematical foundations of Electromagnetism and Gravitation, thus of General Relativity and Gauge Theory, in particular Maxwell and Einstein equations, even if they seem apparently well established.  \\

\noindent
{\bf REFERENCES}  \\

\noindent
[1] BARAKAT, M.: Purity Filtration and the Fine structure of Autonomy, Proceedings of the 19th International Symposium on Mathematical Theory of Networks and Systems, 2010, p 1657-1661.  \\
\noindent
[2] BJORK, J.E.: Analytic D-modules and Applications, Kluwer, 1993.\\
\noindent
[3] BOURBAKI, N.: Alg\`{e}bre homologique, chap. X, Masson, Paris, 1980.  \\
\noindent
[4] BOURBAKI, N.: Alg\`{e}bre commutative, chap. I-IV, Masson, Paris, 1985.  \\
\noindent
[5] CAHEN, M.,GUTT, S.: Invariance des Equations de Maxwell, Bulletin de la Soc. Math\'{e}matique de Belgique, t. XXXIII (1981) 91-97 (See also: M. Cahen et al.(eds), Differential Geometry and Mathematical Physics, Reidel (1983) 27-29 for a translation).  \\
\noindent
[6] COSSERAT, E., COSSERAT, F.: Th\'{e}orie des Corps D\'{e}formables, Hermann, Paris, 1909.\\
\noindent
[7] EISENHART, L.P.: Riemannian Geometry, Princeton University Press, 1926.  \\
\noindent
[8] FOSTER, J., NIGHTINGALE, J.D.:A Short Course in General relativity, Longman, 1979.  \\
\noindent
[9] GR\"{O}BNER, W.: \"{U}ber die Algebraischen Eigenschaften der Integrale von Linearen Differentialgleichungen mit Konstanten Koeffizienten, 
Monatsh. der Math., 47 (1939) 247-284.\\
\noindent
[10] HU, S.-T.: Introduction to Homological Algebra, Holden-Day, 1968.  \\
\noindent
[11] JANET, M.: Sur les Syst\`emes aux d\'eriv\'ees partielles, Journal de Math., 8, 3 (1920) 65-151.\\
\noindent
[12] KALMAN, E.R., YO, Y.C., NARENDA, K.S.: Controllability of Linear Dynamical Systems, Contrib. Diff. Equations, 1, 2 (1963) 189-213.\\
\noindent
[13] KASHIWARA, M.: Algebraic Study of Systems of Partial Differential Equations, M\'emoires de la Soci\'et\'e Math\'ematique de France 63, 1995, 
(Transl. from Japanese of his 1970 Master's Thesis).\\
\noindent
[14] KUNZ, E.: Introduction to Commutative Algebra and Algebraic Geometry, Birkh\"{a}user, 1985.\\
\noindent
[15] MACAULAY, F.S.: The Algebraic Theory of Modular Systems, Cambridge Tracts, vol. 19, Cambridge University Press, London, 1916. Stechert-Hafner Service Agency, New-York, 1964.\\
\noindent
[16] MAISONOBE, P., SABBAH, C.: D-Modules Coh\'erents et Holonomes, Travaux en Cours, 45, Hermann, Paris, 1993.\\
\noindent
[17] NORTHCOTT, D.G.: An Introduction to Homological Algebra, Cambridge University Press, Cambridge, 1966.   \\
\noindent
[18] NORTHCOTT, D.G.: Lessons on Rings, Modules and Multiplicities, Cambridge University Press, Cambridge, 1968.   \\
\noindent
[19] POMMARET, J.-F.: Systems of Partial Differential Equations and Lie Pseudogroups, Gordon and Breach, New York, 1978 (Russian translation: MIR, Moscow, 1983).\\
\noindent
[20] POMMARET, J.-F.: Differential Galois Theory, Gordon and Breach, New York, 1983.\\
\noindent
[21] POMMARET, J.-F.: Lie Pseudogroups and Mechanics, Gordon and Breach, New York, 1988.\\
\noindent
[22] POMMARET, J.-F.: Partial Differential Equations and Group Theory: New Perspectives for Applications, Kluwer, 1994.\\
http://dx.doi.org/10.1007/978-94-017-2539-2  \\
\noindent
[23] POMMARET, J.-F.:Dualit\'e Diff\'erentielle et Applications, C. R. Acad. Sci. Paris, 320, S\'erie I (1995) 1225-1230.\\
\noindent 
[24] POMMARET, J.-F.: Partial Differential Control Theory, Kluwer, 2001 (Zbl 1079.93001).\\
\noindent
[25] POMMARET, J.-F.: Algebraic Analysis of Control Systems Defined by Partial Differential Equations, in Advanced Topics in Control Systems Theory, Lecture Notes in Control and Information Sciences LNCIS 311, Chapter 5, Springer, 2005, 155-223.\\
\noindent
[26] POMMARET, J.-F.: Parametrization of Cosserat Equations, Acta Mechanica, 215 (2010) 43-55.\\
\noindent
[27] POMMARET, J.-F.: Macaulay Inverse Systems Revisited, Journal of Symbolic Computation, 46 (2011) 1049-1069.\\
\noindent
[28] J.-F. POMMARET: Spencer Operator and Applications: From Continuum Mechanics to Mathematical Physics, in "Continuum Mechanics-Progress in Fundamentals and Engineering Applications", Dr. Yong Gan (Ed.), ISBN: 978-953-51-0447--6, InTech, 2012, Available from: \\
http://www.intechopen.com/books/continuum-mechanics-progress-in-fundamentals-and-engineerin-applications/spencer-operator-and-applications-from-continuum-mechanics-to-mathematical-physics  \\
\noindent  
[29] POMMARET, J.-F.: Relative Parametrization of Linear Multidimensional Systems, Multidim. Syst. Sign. Process. (MSSP), Springer, 26 
(2013) 405-437. \\
http://dx.doi.org/10.1007/s11045-013-0265-0  \\
\noindent
[30] POMMARET, J.-F.: The Mathematical Foundations of General Relativity Revisited, Journal of Modern Physics, 4 (2013) 223-239.\\
http://dx.doi.org/10.4236/jmp.2013.48A022   \\
\noindent
[31] POMMARET, J.-F.: The Mathematical Foundations of Gauge Theory Revisited, Journal of Modern Physics, 5 (2014) 157-170.  \\
http://dx.doi.org/10.4236/jmp.2014.55026    \\
\noindent
[32] POMMARET, J.-F.: Deformation Theory of Algebraic and Geometric Structures, Lambert Academic Publisher, (LAP), 
Saarbrucken, Germany, 2016.  \\
http://arxiv.org/abs/1207.1964  \\
\noindent
[33] POMMARET, J.-F.: Airy, Beltrami, Maxwell, Einstein and Lanczos Potentials Revisited, Journal of Modern Physics, 7 (2016) 699-728.  \\
http://dx.doi.org./104236/jmp.2016.77068   \\
\noindent
[34] POMMARET, J.-F.: From Thermodynamics to Gauge Theory: The Virial Theorem Revivited, L. Bailey ed., Gauge Theories and Differential geometry, NOVA publishers, New York, 2016, 1-44. \:
\noindent
[35] POMMARET, J.-F.: Why Gravitational Waves Cannot Exist, Journal of Modern Physics, 8,13 (2017) 2122-2158.  \\
http://dx.doi.org/10.4236/jmp.2017.813130   \\
\noindent
[36] POMMARET, J.-F.: Algebraic Analysis and Mathematical physics, 2017.  \\
http://arxiv.org/abs/1706.04105  \\
\noindent
[37] POMMARET, J.-F.: From Elasticity to Electromagnetism: Beyond the Mirror, \\
http://arxiv.org/abs/1802.02430  \\
\noindent
[38] POMMARET, J.-F.: New Mathematical Methods for Physics, NOVA Science Publisher, New York, 2018.  \\
\noindent
[39] POMMARET, J.-F.: Generating Compatibility Conditions and General Relativity, Journal of Modern Physics, 10, 3 (2019) 371-401.  \\
https://doi.org/10.4236/jmp.2019.103025   \\
\noindent
[40] POMMARET, J.-F., QUADRAT, A.: Localization and parametrization of linear multidimensional control systems, Systems \& Control Letters 37 (1999) 247-260.  \\
\noindent
[41] POMMARET, J.-F., QUADRAT, A.: Algebraic Analysis of Linear Multidimensional  Control Systems, IMA Journal of Mathematical Control and Informations, 16, 1999, 275-297.\\
\noindent
[42] QUADRAT, A.: Analyse Alg\'ebrique des Syst\`emes de Contr\^ole Lin\'eaires Multidimensionnels, Th\`ese de Docteur de l'Ecole Nationale 
des Ponts et Chauss\'ees, 1999 \\
(http://www-sop.inria.fr/cafe/Alban.Quadrat/index.html).\\
\noindent
[43] QUADRAT, A.: Une Introduction \`{a} l'Analyse Alg\'{e}brique Constructive et \`{a} ses Applications, INRIA Research Report 7354, AT-SOP Project, july 2010. Les Cours du CIRM, 1 no. 2: Journ\'{e}es Nationales de Calcul Formel (2010), p281-471 (doi:10.5802/ccirm.11). \\
\noindent
[44] QUADRAT, A.: Grade Filtration of Linear Functional Systems, Acta Applicandae Mathematicae, 127 (2013) 27-86, DOI: 10.1007/s10440-012-9791-2 
(See also  http://hal.inria.fr/inria-00632281/fr/   and  http://pages.saclay.inria.fr/alban.quadrat/PurityFiltration.html ).  \\
\noindent
[45] ROTMAN, J.J.: An Introduction to Homological Algebra, Pure and Applied Mathematics, Academic Press, 1979.\\
\noindent
[46] SCHNEIDERS, J.-P.: ( 1994): An Introduction to D-Modules, Bull. Soc. Roy. Sci. Li\'{e}ge, 63, 223-295.  \\
\noindent
[47] SEILER, W.M.: Involution: The Formal Theory of Differential Equations and its Applications to Computer Algebra, Springer, 2009, 660 pp. (See also doi:10.3842/SIGMA.2009.092 for a recent presentation of involution, in particular sections 3 (p 6 and reference [11], [22]). \\
\noindent
[48] SPENCER, D.C.: Overdetermined Systems of Partial Differential Equations, Bull. Amer. Math. Soc., 75 (1965) 1-114.\\
\noindent
[49] ZERZ, E.: Topics in Multidimensional Linear Systems Theory, Lecture Notes in Control and Information Sciences, LNCIS 256, Springer, 2000.\\

\end{document}